\documentclass[10pt, a4paper, fleqn]{article}
\usepackage{titling}
\title{Unconditional Bound-Preserving and Energy-Dissipating Finite-Volume Schemes\protect\\for the Cahn-Hilliard Equation}

\newcommand{\authorPDF}{Bailo, Carrillo, Kalliadasis, Perez.}
\newcommand{\subjectPDF}{65M08; 35Q92; 35Q35; 35Q70.}
\newcommand{\keywordsPDF}{Cahn-Hilliard equation; diffuse interface theory; gradient flow; finite-volume method; bound preservation; energy dissipation.}
 \usepackage{authblk}

\author[1,2]{Rafael Bailo}
\author[1]{Jos\'e A. Carrillo}
\author[3]{\authorcr Serafim Kalliadasis}
\author[3,4]{Sergio P. Perez}

\affil[1]{Mathematical Institute, University of Oxford}
\affil[2]{Univ. Lille, CNRS, UMR 8524 - Laboratoire Paul Painlev\'e}
\affil[3]{Department of Chemical Engineering, Imperial College London}
\affil[4]{Department of Mathematics, Imperial College London}
\affil[ ]{}
\affil[ ]{\textit{
		bailo@maths.ox.ac.uk,
		carrillo@maths.ox.ac.uk,
	}}
\affil[ ]{\textit{
		s.kalliadasis@imperial.ac.uk,
		sergio.perez15@imperial.ac.uk
	}}
 
\usepackage[left=1.1in, right=1.1in, top=1.1in, bottom=1.1in]{geometry}

\makeatletter
\let\newtitle\@title
\let\newauthor\@author
\let\newdate\@date
\makeatother
\usepackage{fancyhdr}
\usepackage[
	hypertexnames=false,
	pdftex,
	pdftitle={\thetitle. \authorPDF},
	pdfauthor={\authorPDF},
	pdfsubject={\subjectPDF},
	pdfkeywords={\keywordsPDF},
]{hyperref}
\hypersetup{
	colorlinks=true,
	linkcolor=blue,
	citecolor=blue,
	urlcolor=blue,
	linktocpage
}

\usepackage[utf8]{inputenc}
\usepackage[T1]{fontenc}
\usepackage{amsmath}
\usepackage{amsthm}
\usepackage{amsfonts}
\usepackage{amssymb}
\usepackage{graphicx}
\usepackage{mathtools}
\usepackage{ifthen}
\usepackage{dsfont}

\usepackage[UKenglish]{babel}

\usepackage[dvipsnames]{xcolor}

\definecolor{color1}{RGB}{0, 121, 178}
\definecolor{color2}{RGB}{255, 124, 37}
\definecolor{color3}{RGB}{37, 160, 55}
\definecolor{color4}{RGB}{220, 32, 44}
\definecolor{color5}{RGB}{147, 104, 186}
\definecolor{color6}{RGB}{143, 85, 76}
\definecolor{color7}{RGB}{230, 119, 192}
\definecolor{color8}{RGB}{127, 127, 127}
\definecolor{color9}{RGB}{192, 188, 55}
\definecolor{color10}{RGB}{0, 191, 206} 

\newcounter{review}

\newcommand{\ntcreview}[3]{\refstepcounter{review}

	{\color{#2}{\textbf{[#1]}: #3}}}

\newcommand{\creview}[3]{\ntcreview{#1}{#2}{#3}
	\addcontentsline{tor}{subsection}{\thereview~\textbf{[#1]}:~#3
	}}

\newcommand{\review}[2]{\creview{#1}{blue}{#2}}

\makeatletter

\newcommand\listreviewname{List of Reviews}
\newcommand\listofreviews{\section*{\listreviewname}\@starttoc{tor}}
\makeatother

\usepackage[all]{nowidow}

\newcommand{\subjectclassification}[1]{

	{\small\textbf{\textit{AMS Subject Classification --- }} #1}

}
\newcommand{\keywords}[1]{

	{\small\textbf{\textit{Keywords --- }} #1}

}

\renewcommand\lll\MoveEqLeft

\usepackage{tikz}
\usetikzlibrary{patterns}
\usetikzlibrary{decorations.pathreplacing}
\usetikzlibrary{calc}
\usetikzlibrary{fadings}
\usetikzlibrary{shapes}
\usetikzlibrary{plotmarks}
\usetikzlibrary{arrows, decorations.markings}
\tikzset{thicker line small arrows m/.style args={#1in#2}{
			draw=#2,
			solid,
			line width=#1,
			shorten >=1mm,
			decoration={
					markings,
					mark=at position 1.0 with {\arrow[fill=#2,thin]{triangle 90}}
				},
			postaction={decorate}
		}}

\usepackage{pgfplots}
\pgfplotsset{compat=1.16}

\usepackage{float}
\usepackage{subcaption}
\usepackage[section]{placeins}
\usepackage{rotating}

\usepackage{array}
\newcolumntype{L}[1]{>{\raggedright\let\newline\\\arraybackslash\hspace{0pt}}m{#1}}
\newcolumntype{C}[1]{>{\centering\let\newline\\\arraybackslash\hspace{0pt}}m{#1}}
\newcolumntype{R}[1]{>{\raggedleft\let\newline\\\arraybackslash\hspace{0pt}}m{#1}}
\usepackage{booktabs}
\usepackage{tabularx}
\usepackage{multirow}

\usepackage{titling}

\usepackage[algoruled]{algorithm2e}

\usepackage[capitalise]{cleveref}

 \usepackage{accents}

\newcommand\term\emph

\numberwithin{equation}{section}

\usepackage{autobreak}

\usepackage{enumerate}

\makeatletter
\def\@maketitle{\newpage
	\begin{center}\let \footnote \thanks
		{\LARGE\bfseries \@title \par}\vskip 2.5em{\large
				\lineskip .5em\begin{tabular}[t]{c}\@author
				\end{tabular}\par}\vskip 1em{\large \@date}\end{center}\par
	\vskip 1.5em}
\makeatother



\usepackage{amsthm}
\theoremstyle{plain}
\newtheorem{theorem}{Theorem}[section]
\newtheorem{lemma}[theorem]{Lemma}

\theoremstyle{remark}
\newtheorem{remark}[theorem]{\bf Remark}

\usepackage{esint}

\def\XXint#1#2#3{{\setbox0=\hbox{$#1{#2#3}{\int}$ }
			\vcenter{\hbox{$#2#3$ }}\kern-.6\wd0}}

\DeclarePairedDelimiter{\prt}{(}{)}
\DeclarePairedDelimiter{\brk}{[}{]}

\DeclarePairedDelimiter{\abs}{|}{|}
\DeclarePairedDelimiter{\norm}{\|}{\|}
\DeclarePairedDelimiter{\set}{\{}{\}}

\DeclarePairedDelimiter{\inn}{\langle}{\rangle}

\usepackage{suffix}

\newcommand{\inner}[2]{\inn{#1,#2}}
\WithSuffix\newcommand\inner*[2]{\inn*{#1,#2}}

\DeclarePairedDelimiter{\positive}{(}{)^{+}}
\DeclarePairedDelimiter{\negative}{(}{)^{-}}

\newcommand\pos\positive
\renewcommand\neg\negative

\WithSuffix\newcommand\pos*{\positive*}
\WithSuffix\newcommand\neg*{\negative*}

\renewcommand{\L}[1]{{L^{#1}}}

\newcommand{\pnorm}[2]{\norm{#2}_{\L{#1}}}
\WithSuffix\newcommand\pnorm*[2]{\norm*{#2}_{\L{#1}}}

\newcommand{\psnorm}[3]{\norm{#3}_{\L{#1}(#2)}}
\WithSuffix\newcommand\psnorm*[3]{\norm*{#3}_{\L{#1}(#2)}}

\newcommand{\pnormp}[2]{\pnorm{#1}{#2}^{#1}}
\WithSuffix\newcommand\pnormp*[2]{\pnorm*{#1}{#2}^{#1}}

\newcommand{\psnormp}[3]{\psnorm{#1}{#2}{#3}^{#1}}
\WithSuffix\newcommand\psnormp*[3]{\psnorm*{#1}{#2}{#3}^{#1}}

\newcommand\svec\vec

\renewcommand{\vec}{\mathbf}
\renewcommand{\svec}{\boldsymbol}

\renewcommand{\d}{\mathrm{d}}

\newcommand{\dd}{\mathop{}\!\d}

\newcommand{\der}[2]{\frac{\d #1}{\d #2}}
\newcommand{\pder}[2]{\frac{\partial #1}{\partial #2}}

\newcommand{\vder}[2]{\frac{\delta #1}{\delta #2}}

\newcommand{\ds}{\dd s}
\newcommand{\dx}{\dd x}
\newcommand{\dy}{\dd y}

\newcommand{\grad}{\nabla}
\renewcommand{\div}{\nabla\cdot}

\newcommand{\laplacian}{\Delta}
\newcommand{\laplace}{\laplacian}

\newcommand{\Dt}{\Delta t}
\newcommand{\Dx}{\Delta x}
\newcommand{\Dy}{\Delta y}

\newcommand{\nhalf}{1/2}

\renewcommand{\i}{_{i}}
\newcommand{\ip}{_{i+1}}

\newcommand{\ih}{_{i+\nhalf}}
\newcommand{\imh}{_{i-\nhalf}}

\renewcommand{\j}{_{j}}

\newcommand{\jh}{_{j+\nhalf}}
\newcommand{\jmh}{_{j-\nhalf}}

\renewcommand{\ij}{_{i,\,j}}
\newcommand{\ipj}{_{i+1,\,j}}
\newcommand{\imj}{_{i-1,\,j}}
\newcommand{\ihj}{_{i+\nhalf,\,j}}
\newcommand{\imhj}{_{i-\nhalf,\,j}}

\newcommand{\ijp}{_{i,\,j+1}}
\newcommand{\ijm}{_{i,\,j-1}}
\newcommand{\ijh}{_{i,\,j+\nhalf}}
\newcommand{\ijmh}{_{i,\,j-\nhalf}}

\newcommand{\n}{^{n}}
\newcommand{\np}{^{n+1}}

\newcommand{\ppr}{r}
\newcommand{\pprm}{r-1}

\newcommand{\nr}{^{n,\,\ppr}}

\newcommand{\nrm}{^{n,\,\pprm}}

\newcommand{\Wr}{^{W,\,\ppr}}

\newcommand{\nh}{^{n+\nhalf}}

\newlength{\dhatheight}

\ifthenelse{\isundefined{\Wr}}{
	\newcommand{\Wr}{^{W,\,\ppr}}
}{
	\renewcommand{\Wr}{^{W,\,\ppr}}
}

 \usepackage{todonotes}

\DeclareCaptionFormat{cont}{#1 (cont.)#2#3\par}

\newif\ifskiptable

\pgfplotsset{colormap={hsv}{
			hsb(0.00cm)=(0.00,0,0.95);
			hsb(0.05cm)=(0.05,1,1);
			hsb(0.10cm)=(0.10,1,1);
			hsb(0.15cm)=(0.15,1,1);
			hsb(0.20cm)=(0.20,1,1);
			hsb(0.25cm)=(0.25,1,1);
			hsb(0.30cm)=(0.30,1,1);
			hsb(0.35cm)=(0.35,1,1);
			hsb(0.40cm)=(0.40,1,1);
			hsb(0.45cm)=(0.45,1,1);
			hsb(0.50cm)=(0.50,1,1);
			hsb(0.55cm)=(0.55,1,1);
			hsb(0.60cm)=(0.60,1,1);
			hsb(0.65cm)=(0.65,1,1);
			hsb(0.70cm)=(0.70,1,1);
			hsb(0.75cm)=(0.75,1,1);
			hsb(0.80cm)=(0.80,1,1);
			hsb(0.85cm)=(0.85,1,1);
			hsb(0.90cm)=(0.90,1,1);
			hsb(0.95cm)=(0.95,1,1);
			hsb(1.00cm)=(1.00,1,1);
		}
}

\pgfplotsset{colormap={hsvSoft}{
			hsb(0.00cm)=(0.00,0,0.95);
			hsb(0.05cm)=(0.05,1,1);
			hsb(0.10cm)=(0.10,1,1);
			hsb(0.15cm)=(0.15,1,1);
			hsb(0.20cm)=(0.20,1,1);
			hsb(0.25cm)=(0.25,1,1);
			hsb(0.30cm)=(0.30,1,1);
			hsb(0.35cm)=(0.35,1,1);
			hsb(0.40cm)=(0.40,1,1);
			hsb(0.45cm)=(0.45,1,1);
			hsb(0.50cm)=(0.50,1,1);
			hsb(0.55cm)=(0.55,1,1);
			hsb(0.60cm)=(0.60,1,1);
			hsb(0.65cm)=(0.65,1,1);
			hsb(0.70cm)=(0.70,1,1);
			hsb(0.75cm)=(0.75,1,1);
			hsb(0.80cm)=(0.80,1,1);
			hsb(0.85cm)=(0.85,1,1);
			hsb(0.90cm)=(0.90,1,1);
			hsb(0.95cm)=(0.95,1,1);
			hsb(1.00cm)=(0.00,0,0.95);
		}
}

\pgfplotsset{colormap={viridisFull}{
			rgb=(0.26700401, 0.00487433, 0.32941519)
			rgb=(0.26851048, 0.00960483, 0.33542652)
			rgb=(0.26994384, 0.01462494, 0.34137895)
			rgb=(0.27130489, 0.01994186, 0.34726862)
			rgb=(0.27259384, 0.02556309, 0.35309303)
			rgb=(0.27380934, 0.03149748, 0.35885256)
			rgb=(0.27495242, 0.03775181, 0.36454323)
			rgb=(0.27602238, 0.04416723, 0.37016418)
			rgb=(0.2770184 , 0.05034437, 0.37571452)
			rgb=(0.27794143, 0.05632444, 0.38119074)
			rgb=(0.27879067, 0.06214536, 0.38659204)
			rgb=(0.2795655 , 0.06783587, 0.39191723)
			rgb=(0.28026658, 0.07341724, 0.39716349)
			rgb=(0.28089358, 0.07890703, 0.40232944)
			rgb=(0.28144581, 0.0843197 , 0.40741404)
			rgb=(0.28192358, 0.08966622, 0.41241521)
			rgb=(0.28232739, 0.09495545, 0.41733086)
			rgb=(0.28265633, 0.10019576, 0.42216032)
			rgb=(0.28291049, 0.10539345, 0.42690202)
			rgb=(0.28309095, 0.11055307, 0.43155375)
			rgb=(0.28319704, 0.11567966, 0.43611482)
			rgb=(0.28322882, 0.12077701, 0.44058404)
			rgb=(0.28318684, 0.12584799, 0.44496 )
			rgb=(0.283072 , 0.13089477, 0.44924127)
			rgb=(0.28288389, 0.13592005, 0.45342734)
			rgb=(0.28262297, 0.14092556, 0.45751726)
			rgb=(0.28229037, 0.14591233, 0.46150995)
			rgb=(0.28188676, 0.15088147, 0.46540474)
			rgb=(0.28141228, 0.15583425, 0.46920128)
			rgb=(0.28086773, 0.16077132, 0.47289909)
			rgb=(0.28025468, 0.16569272, 0.47649762)
			rgb=(0.27957399, 0.17059884, 0.47999675)
			rgb=(0.27882618, 0.1754902 , 0.48339654)
			rgb=(0.27801236, 0.18036684, 0.48669702)
			rgb=(0.27713437, 0.18522836, 0.48989831)
			rgb=(0.27619376, 0.19007447, 0.49300074)
			rgb=(0.27519116, 0.1949054 , 0.49600488)
			rgb=(0.27412802, 0.19972086, 0.49891131)
			rgb=(0.27300596, 0.20452049, 0.50172076)
			rgb=(0.27182812, 0.20930306, 0.50443413)
			rgb=(0.27059473, 0.21406899, 0.50705243)
			rgb=(0.26930756, 0.21881782, 0.50957678)
			rgb=(0.26796846, 0.22354911, 0.5120084 )
			rgb=(0.26657984, 0.2282621 , 0.5143487 )
			rgb=(0.2651445 , 0.23295593, 0.5165993 )
			rgb=(0.2636632 , 0.23763078, 0.51876163)
			rgb=(0.26213801, 0.24228619, 0.52083736)
			rgb=(0.26057103, 0.2469217 , 0.52282822)
			rgb=(0.25896451, 0.25153685, 0.52473609)
			rgb=(0.25732244, 0.2561304 , 0.52656332)
			rgb=(0.25564519, 0.26070284, 0.52831152)
			rgb=(0.25393498, 0.26525384, 0.52998273)
			rgb=(0.25219404, 0.26978306, 0.53157905)
			rgb=(0.25042462, 0.27429024, 0.53310261)
			rgb=(0.24862899, 0.27877509, 0.53455561)
			rgb=(0.2468114 , 0.28323662, 0.53594093)
			rgb=(0.24497208, 0.28767547, 0.53726018)
			rgb=(0.24311324, 0.29209154, 0.53851561)
			rgb=(0.24123708, 0.29648471, 0.53970946)
			rgb=(0.23934575, 0.30085494, 0.54084398)
			rgb=(0.23744138, 0.30520222, 0.5419214 )
			rgb=(0.23552606, 0.30952657, 0.54294396)
			rgb=(0.23360277, 0.31382773, 0.54391424)
			rgb=(0.2316735 , 0.3181058 , 0.54483444)
			rgb=(0.22973926, 0.32236127, 0.54570633)
			rgb=(0.22780192, 0.32659432, 0.546532 )
			rgb=(0.2258633 , 0.33080515, 0.54731353)
			rgb=(0.22392515, 0.334994 , 0.54805291)
			rgb=(0.22198915, 0.33916114, 0.54875211)
			rgb=(0.22005691, 0.34330688, 0.54941304)
			rgb=(0.21812995, 0.34743154, 0.55003755)
			rgb=(0.21620971, 0.35153548, 0.55062743)
			rgb=(0.21429757, 0.35561907, 0.5511844 )
			rgb=(0.21239477, 0.35968273, 0.55171011)
			rgb=(0.2105031 , 0.36372671, 0.55220646)
			rgb=(0.20862342, 0.36775151, 0.55267486)
			rgb=(0.20675628, 0.37175775, 0.55311653)
			rgb=(0.20490257, 0.37574589, 0.55353282)
			rgb=(0.20306309, 0.37971644, 0.55392505)
			rgb=(0.20123854, 0.38366989, 0.55429441)
			rgb=(0.1994295 , 0.38760678, 0.55464205)
			rgb=(0.1976365 , 0.39152762, 0.55496905)
			rgb=(0.19585993, 0.39543297, 0.55527637)
			rgb=(0.19410009, 0.39932336, 0.55556494)
			rgb=(0.19235719, 0.40319934, 0.55583559)
			rgb=(0.19063135, 0.40706148, 0.55608907)
			rgb=(0.18892259, 0.41091033, 0.55632606)
			rgb=(0.18723083, 0.41474645, 0.55654717)
			rgb=(0.18555593, 0.4185704 , 0.55675292)
			rgb=(0.18389763, 0.42238275, 0.55694377)
			rgb=(0.18225561, 0.42618405, 0.5571201 )
			rgb=(0.18062949, 0.42997486, 0.55728221)
			rgb=(0.17901879, 0.43375572, 0.55743035)
			rgb=(0.17742298, 0.4375272 , 0.55756466)
			rgb=(0.17584148, 0.44128981, 0.55768526)
			rgb=(0.17427363, 0.4450441 , 0.55779216)
			rgb=(0.17271876, 0.4487906 , 0.55788532)
			rgb=(0.17117615, 0.4525298 , 0.55796464)
			rgb=(0.16964573, 0.45626209, 0.55803034)
			rgb=(0.16812641, 0.45998802, 0.55808199)
			rgb=(0.1666171 , 0.46370813, 0.55811913)
			rgb=(0.16511703, 0.4674229 , 0.55814141)
			rgb=(0.16362543, 0.47113278, 0.55814842)
			rgb=(0.16214155, 0.47483821, 0.55813967)
			rgb=(0.16066467, 0.47853961, 0.55811466)
			rgb=(0.15919413, 0.4822374 , 0.5580728 )
			rgb=(0.15772933, 0.48593197, 0.55801347)
			rgb=(0.15626973, 0.4896237 , 0.557936 )
			rgb=(0.15481488, 0.49331293, 0.55783967)
			rgb=(0.15336445, 0.49700003, 0.55772371)
			rgb=(0.1519182 , 0.50068529, 0.55758733)
			rgb=(0.15047605, 0.50436904, 0.55742968)
			rgb=(0.14903918, 0.50805136, 0.5572505 )
			rgb=(0.14760731, 0.51173263, 0.55704861)
			rgb=(0.14618026, 0.51541316, 0.55682271)
			rgb=(0.14475863, 0.51909319, 0.55657181)
			rgb=(0.14334327, 0.52277292, 0.55629491)
			rgb=(0.14193527, 0.52645254, 0.55599097)
			rgb=(0.14053599, 0.53013219, 0.55565893)
			rgb=(0.13914708, 0.53381201, 0.55529773)
			rgb=(0.13777048, 0.53749213, 0.55490625)
			rgb=(0.1364085 , 0.54117264, 0.55448339)
			rgb=(0.13506561, 0.54485335, 0.55402906)
			rgb=(0.13374299, 0.54853458, 0.55354108)
			rgb=(0.13244401, 0.55221637, 0.55301828)
			rgb=(0.13117249, 0.55589872, 0.55245948)
			rgb=(0.1299327 , 0.55958162, 0.55186354)
			rgb=(0.12872938, 0.56326503, 0.55122927)
			rgb=(0.12756771, 0.56694891, 0.55055551)
			rgb=(0.12645338, 0.57063316, 0.5498411 )
			rgb=(0.12539383, 0.57431754, 0.54908564)
			rgb=(0.12439474, 0.57800205, 0.5482874 )
			rgb=(0.12346281, 0.58168661, 0.54744498)
			rgb=(0.12260562, 0.58537105, 0.54655722)
			rgb=(0.12183122, 0.58905521, 0.54562298)
			rgb=(0.12114807, 0.59273889, 0.54464114)
			rgb=(0.12056501, 0.59642187, 0.54361058)
			rgb=(0.12009154, 0.60010387, 0.54253043)
			rgb=(0.11973756, 0.60378459, 0.54139999)
			rgb=(0.11951163, 0.60746388, 0.54021751)
			rgb=(0.11942341, 0.61114146, 0.53898192)
			rgb=(0.11948255, 0.61481702, 0.53769219)
			rgb=(0.11969858, 0.61849025, 0.53634733)
			rgb=(0.12008079, 0.62216081, 0.53494633)
			rgb=(0.12063824, 0.62582833, 0.53348834)
			rgb=(0.12137972, 0.62949242, 0.53197275)
			rgb=(0.12231244, 0.63315277, 0.53039808)
			rgb=(0.12344358, 0.63680899, 0.52876343)
			rgb=(0.12477953, 0.64046069, 0.52706792)
			rgb=(0.12632581, 0.64410744, 0.52531069)
			rgb=(0.12808703, 0.64774881, 0.52349092)
			rgb=(0.13006688, 0.65138436, 0.52160791)
			rgb=(0.13226797, 0.65501363, 0.51966086)
			rgb=(0.13469183, 0.65863619, 0.5176488 )
			rgb=(0.13733921, 0.66225157, 0.51557101)
			rgb=(0.14020991, 0.66585927, 0.5134268 )
			rgb=(0.14330291, 0.66945881, 0.51121549)
			rgb=(0.1466164 , 0.67304968, 0.50893644)
			rgb=(0.15014782, 0.67663139, 0.5065889 )
			rgb=(0.15389405, 0.68020343, 0.50417217)
			rgb=(0.15785146, 0.68376525, 0.50168574)
			rgb=(0.16201598, 0.68731632, 0.49912906)
			rgb=(0.1663832 , 0.69085611, 0.49650163)
			rgb=(0.1709484 , 0.69438405, 0.49380294)
			rgb=(0.17570671, 0.6978996 , 0.49103252)
			rgb=(0.18065314, 0.70140222, 0.48818938)
			rgb=(0.18578266, 0.70489133, 0.48527326)
			rgb=(0.19109018, 0.70836635, 0.48228395)
			rgb=(0.19657063, 0.71182668, 0.47922108)
			rgb=(0.20221902, 0.71527175, 0.47608431)
			rgb=(0.20803045, 0.71870095, 0.4728733 )
			rgb=(0.21400015, 0.72211371, 0.46958774)
			rgb=(0.22012381, 0.72550945, 0.46622638)
			rgb=(0.2263969 , 0.72888753, 0.46278934)
			rgb=(0.23281498, 0.73224735, 0.45927675)
			rgb=(0.2393739 , 0.73558828, 0.45568838)
			rgb=(0.24606968, 0.73890972, 0.45202405)
			rgb=(0.25289851, 0.74221104, 0.44828355)
			rgb=(0.25985676, 0.74549162, 0.44446673)
			rgb=(0.26694127, 0.74875084, 0.44057284)
			rgb=(0.27414922, 0.75198807, 0.4366009 )
			rgb=(0.28147681, 0.75520266, 0.43255207)
			rgb=(0.28892102, 0.75839399, 0.42842626)
			rgb=(0.29647899, 0.76156142, 0.42422341)
			rgb=(0.30414796, 0.76470433, 0.41994346)
			rgb=(0.31192534, 0.76782207, 0.41558638)
			rgb=(0.3198086 , 0.77091403, 0.41115215)
			rgb=(0.3277958 , 0.77397953, 0.40664011)
			rgb=(0.33588539, 0.7770179 , 0.40204917)
			rgb=(0.34407411, 0.78002855, 0.39738103)
			rgb=(0.35235985, 0.78301086, 0.39263579)
			rgb=(0.36074053, 0.78596419, 0.38781353)
			rgb=(0.3692142 , 0.78888793, 0.38291438)
			rgb=(0.37777892, 0.79178146, 0.3779385 )
			rgb=(0.38643282, 0.79464415, 0.37288606)
			rgb=(0.39517408, 0.79747541, 0.36775726)
			rgb=(0.40400101, 0.80027461, 0.36255223)
			rgb=(0.4129135 , 0.80304099, 0.35726893)
			rgb=(0.42190813, 0.80577412, 0.35191009)
			rgb=(0.43098317, 0.80847343, 0.34647607)
			rgb=(0.44013691, 0.81113836, 0.3409673 )
			rgb=(0.44936763, 0.81376835, 0.33538426)
			rgb=(0.45867362, 0.81636288, 0.32972749)
			rgb=(0.46805314, 0.81892143, 0.32399761)
			rgb=(0.47750446, 0.82144351, 0.31819529)
			rgb=(0.4870258 , 0.82392862, 0.31232133)
			rgb=(0.49661536, 0.82637633, 0.30637661)
			rgb=(0.5062713 , 0.82878621, 0.30036211)
			rgb=(0.51599182, 0.83115784, 0.29427888)
			rgb=(0.52577622, 0.83349064, 0.2881265 )
			rgb=(0.5356211 , 0.83578452, 0.28190832)
			rgb=(0.5455244 , 0.83803918, 0.27562602)
			rgb=(0.55548397, 0.84025437, 0.26928147)
			rgb=(0.5654976 , 0.8424299 , 0.26287683)
			rgb=(0.57556297, 0.84456561, 0.25641457)
			rgb=(0.58567772, 0.84666139, 0.24989748)
			rgb=(0.59583934, 0.84871722, 0.24332878)
			rgb=(0.60604528, 0.8507331 , 0.23671214)
			rgb=(0.61629283, 0.85270912, 0.23005179)
			rgb=(0.62657923, 0.85464543, 0.22335258)
			rgb=(0.63690157, 0.85654226, 0.21662012)
			rgb=(0.64725685, 0.85839991, 0.20986086)
			rgb=(0.65764197, 0.86021878, 0.20308229)
			rgb=(0.66805369, 0.86199932, 0.19629307)
			rgb=(0.67848868, 0.86374211, 0.18950326)
			rgb=(0.68894351, 0.86544779, 0.18272455)
			rgb=(0.69941463, 0.86711711, 0.17597055)
			rgb=(0.70989842, 0.86875092, 0.16925712)
			rgb=(0.72039115, 0.87035015, 0.16260273)
			rgb=(0.73088902, 0.87191584, 0.15602894)
			rgb=(0.74138803, 0.87344918, 0.14956101)
			rgb=(0.75188414, 0.87495143, 0.14322828)
			rgb=(0.76237342, 0.87642392, 0.13706449)
			rgb=(0.77285183, 0.87786808, 0.13110864)
			rgb=(0.78331535, 0.87928545, 0.12540538)
			rgb=(0.79375994, 0.88067763, 0.12000532)
			rgb=(0.80418159, 0.88204632, 0.11496505)
			rgb=(0.81457634, 0.88339329, 0.11034678)
			rgb=(0.82494028, 0.88472036, 0.10621724)
			rgb=(0.83526959, 0.88602943, 0.1026459 )
			rgb=(0.84556056, 0.88732243, 0.09970219)
			rgb=(0.8558096 , 0.88860134, 0.09745186)
			rgb=(0.86601325, 0.88986815, 0.09595277)
			rgb=(0.87616824, 0.89112487, 0.09525046)
			rgb=(0.88627146, 0.89237353, 0.09537439)
			rgb=(0.89632002, 0.89361614, 0.09633538)
			rgb=(0.90631121, 0.89485467, 0.09812496)
			rgb=(0.91624212, 0.89609127, 0.1007168 )
			rgb=(0.92610579, 0.89732977, 0.10407067)
			rgb=(0.93590444, 0.8985704 , 0.10813094)
			rgb=(0.94563626, 0.899815 , 0.11283773)
			rgb=(0.95529972, 0.90106534, 0.11812832)
			rgb=(0.96489353, 0.90232311, 0.12394051)
			rgb=(0.97441665, 0.90358991, 0.13021494)
			rgb=(0.98386829, 0.90486726, 0.13689671)
			rgb=(0.99324789, 0.90615657, 0.1439362 )
		}
}

\pgfplotsset{colormap={viridisSoft}{
			rgb255=(242, 242, 242);
rgb=(0.28026,0.1657,0.4765);
			rgb=(0.26366,0.23763,0.51877);
			rgb=(0.23744,0.3052,0.54192);
			rgb=(0.20862,0.36775,0.55267);
			rgb=(0.18225,0.42618,0.55711);
			rgb=(0.1592,0.48224,0.55807);
			rgb=(0.13777,0.53749,0.5549);
			rgb=(0.12115,0.59274,0.54465);
			rgb=(0.12808,0.64775,0.5235);
			rgb=(0.18065,0.7014,0.48819);
			rgb=(0.27415,0.75198,0.4366);
			rgb=(0.39517,0.79747,0.36775);
			rgb=(0.53561,0.83578,0.2819);
			rgb=(0.68895,0.86545,0.18272);
			rgb=(0.84557,0.88733,0.0997);
			rgb=(0.99324,0.90616,0.14394)
		}
}

\pgfplotsset{colormap={cellRed}{
			rgb255=(242.0,242.0,242.0);
			rgb255=(241.63157894736844,234.47368421052633,234.47368421052633);
			rgb255=(241.26315789473685,226.94736842105266,226.94736842105266);
			rgb255=(240.89473684210526,219.42105263157893,219.42105263157893);
			rgb255=(240.5263157894737,211.89473684210526,211.89473684210526);
			rgb255=(240.1578947368421,204.3684210526316,204.3684210526316);
			rgb255=(239.78947368421052,196.84210526315792,196.84210526315792);
			rgb255=(239.42105263157896,189.31578947368422,189.31578947368422);
			rgb255=(239.05263157894737,181.78947368421052,181.78947368421052);
			rgb255=(238.6842105263158,174.26315789473688,174.26315789473688);
			rgb255=(238.31578947368422,166.73684210526315,166.73684210526315);
			rgb255=(237.94736842105263,159.21052631578948,159.21052631578948);
			rgb255=(237.57894736842104,151.68421052631578,151.68421052631578);
			rgb255=(237.21052631578948,144.1578947368421,144.1578947368421);
			rgb255=(236.84210526315792,136.63157894736844,136.63157894736844);
			rgb255=(236.47368421052633,129.10526315789474,129.10526315789474);
			rgb255=(236.10526315789474,121.57894736842107,121.57894736842107);
			rgb255=(235.73684210526318,114.05263157894737,114.05263157894737);
			rgb255=(235.3684210526316,106.52631578947368,106.52631578947368);
			rgb255=(235.0,99.0,99.0);
		}
}

\pgfplotsset{colormap={cellGreen}{
			rgb255=(242.0,242.0,242.0);
			rgb255=(236.21052631578948,239.5263157894737,234.26315789473685);
			rgb255=(230.42105263157896,237.05263157894737,226.5263157894737);
			rgb255=(224.6315789473684,234.57894736842104,218.78947368421052);
			rgb255=(218.8421052631579,232.10526315789474,211.05263157894737);
			rgb255=(213.05263157894737,229.63157894736844,203.31578947368422);
			rgb255=(207.26315789473685,227.1578947368421,195.57894736842107);
			rgb255=(201.4736842105263,224.68421052631578,187.8421052631579);
			rgb255=(195.68421052631578,222.21052631578948,180.10526315789474);
			rgb255=(189.8947368421053,219.73684210526318,172.36842105263162);
			rgb255=(184.10526315789474,217.26315789473682,164.63157894736844);
			rgb255=(178.31578947368422,214.78947368421052,156.89473684210526);
			rgb255=(172.5263157894737,212.31578947368422,149.1578947368421);
			rgb255=(166.73684210526318,209.84210526315792,141.42105263157896);
			rgb255=(160.94736842105263,207.3684210526316,133.6842105263158);
			rgb255=(155.1578947368421,204.89473684210526,125.94736842105263);
			rgb255=(149.3684210526316,202.42105263157893,118.21052631578948);
			rgb255=(143.57894736842104,199.94736842105266,110.47368421052632);
			rgb255=(137.78947368421052,197.47368421052633,102.73684210526316);
			rgb255=(132.0,195.0,95.0);
		}
}

\pgfplotsset{colormap={cellRedSquared}{
			rgb255=(242.0,242.0,242.0);
			rgb255=(241.28254847645428,227.34349030470915,227.34349030470915);
			rgb255=(240.60387811634348,213.47922437673128,213.47922437673128);
			rgb255=(239.9639889196676,200.40720221606648,200.40720221606648);
			rgb255=(239.36288088642658,188.1274238227147,188.1274238227147);
			rgb255=(238.8005540166205,176.63988919667594,176.63988919667594);
			rgb255=(238.2770083102493,165.94459833795014,165.94459833795014);
			rgb255=(237.79224376731304,156.04155124653738,156.04155124653738);
			rgb255=(237.34626038781164,146.93074792243766,146.93074792243766);
			rgb255=(236.93905817174516,138.61218836565098,138.61218836565098);
			rgb255=(236.57063711911357,131.0858725761773,131.0858725761773);
			rgb255=(236.2409972299169,124.35180055401662,124.35180055401662);
			rgb255=(235.95013850415512,118.40997229916897,118.40997229916897);
			rgb255=(235.69806094182823,113.26038781163435,113.26038781163435);
			rgb255=(235.4847645429363,108.90304709141274,108.90304709141274);
			rgb255=(235.3102493074792,105.33795013850416,105.33795013850416);
			rgb255=(235.17451523545705,102.56509695290858,102.56509695290858);
			rgb255=(235.0775623268698,100.58448753462605,100.58448753462605);
			rgb255=(235.01939058171746,99.3961218836565,99.3961218836565);
			rgb255=(235.0,99.0,99.0);
		}
}
\pgfplotsset{colormap={cellGreenSquared}{
			rgb255=(242.0,242.0,242.0);
			rgb255=(230.7257617728532,237.18282548476455,226.93351800554018);
			rgb255=(220.06094182825484,232.62603878116343,212.6814404432133);
			rgb255=(210.00554016620498,228.32963988919667,199.2437673130194);
			rgb255=(200.5595567867036,224.29362880886427,186.62049861495845);
			rgb255=(191.7229916897507,220.5180055401662,174.8116343490305);
			rgb255=(183.49584487534625,217.0027700831025,163.81717451523545);
			rgb255=(175.87811634349032,213.74792243767314,153.63711911357342);
			rgb255=(168.86980609418282,210.75346260387812,144.27146814404432);
			rgb255=(162.47091412742384,208.01939058171746,135.72022160664818);
			rgb255=(156.68144044321332,205.54570637119116,127.98337950138506);
			rgb255=(151.50138504155126,203.33240997229916,121.06094182825484);
			rgb255=(146.9307479224377,201.37950138504155,114.95290858725764);
			rgb255=(142.96952908587255,199.68698060941827,109.65927977839334);
			rgb255=(139.61772853185596,198.25484764542935,105.18005540166205);
			rgb255=(136.8753462603878,197.0831024930748,101.51523545706371);
			rgb255=(134.74238227146813,196.17174515235456,98.66481994459834);
			rgb255=(133.21883656509695,195.5207756232687,96.62880886426592);
			rgb255=(132.30470914127426,195.13019390581718,95.40720221606648);
			rgb255=(132.0,195.0,95.0);
		}
}

\usepgfplotslibrary{patchplots}

\pgfplotsset{every axis/.append style={
			grid=both,
			grid style={white, line width=.1pt},
			major grid style={white, line width=1.5pt},
			axis background/.style={fill=gray!10},
			axis line style={draw=none},
			tick style={draw=none},
			xlabel = $x$,
line width=1pt,
legend style={
					line width = 1pt,
					draw=none,
					/tikz/every even column/.append style={column sep=0.5cm}
				},
		}}

\usepgfplotslibrary{groupplots}

\usepackage{calc}

\definecolor{gg0}{HTML}{E24A33}
\definecolor{gg1}{HTML}{348ABD}
\definecolor{gg2}{HTML}{988ED5}
\definecolor{gg3}{HTML}{777777}
\definecolor{gg4}{HTML}{FBC15E}
\definecolor{gg5}{HTML}{8EBA42}
\definecolor{gg6}{HTML}{FFB5B8}

\pgfplotsset{
	/pgfplots/colormap={bright}{rgb255=(0,0,0) rgb255=(78,3,100) rgb255=(2,74,255)
			rgb255=(255,21,181) rgb255=(255,113,26) rgb255=(147,213,114) rgb255=(230,255,0)
			rgb255=(255,255,255)}
}

\newcommand{\eps}{\varepsilon}

\usepackage[nolist]{acronym}

\newcommand{\F}{\mathcal{F}}

\newcommand{\GL}{_{\text{GL}}}
\newcommand{\FH}{_{\text{log}}}

\newcommand{\nhc}{^{n+\nhalf,\,c}}
\newcommand{\nhcm}{^{n+\nhalf,\,c-1}}

\newcommand{\nnr}{^{r}}
\newcommand{\nnrm}{^{r-1}}

\newcommand{\nnhc}{^{c}}
\newcommand{\nnhcm}{^{c-1}}

\newcommand{\addappendix}{\section*{\appendixname}\addcontentsline{toc}{section}{\appendixname}\counterwithin*{figure}{section}\stepcounter{section}\renewcommand{\thesection}{A}\renewcommand{\thefigure}{\thesection.\arabic{figure}}} 

\allowdisplaybreaks

\renewcommand{\review}[2]{}
\renewcommand{\creview}[3]{}
\renewcommand{\ntcreview}[3]{}

\renewcommand{\tableofcontents}{}
\renewcommand{\listofreviews}{}

\expandafter\def\csname ver@etex.sty\endcsname{3000/12/31}

\usepackage{autonum}
\makeatletter
\patchcmd{\autonum@saveEnvironmentSubcommands}
{(0,0)\begin}
{(0,0)\hfuzz=\maxdimen\begin}
{}{}
\makeatother

\AtBeginDocument{

}

\makeatletter
\autonum@generatePatchedReferenceCSL{ref}
\autonum@generatePatchedReferenceCSL{eqref}
\autonum@generatePatchedReferenceCSL{Cref}
\makeatother

\definecolor{revisionColourOne}{RGB}{180,0,0}
\definecolor{revisionColourTwo}{RGB}{0,0,180}

\usepackage{setspace}
\begin{document}
\begin{singlespace}\maketitle\end{singlespace}
\begin{abstract}
	We propose finite-volume schemes for the Cahn-Hilliard equation which unconditionally and discretely preserve the boundedness of the phase field and the dissipation of the free energy. Our numerical framework is applicable to a variety of free-energy potentials, including Ginzburg-Landau and Flory-Huggins, to general wetting boundary conditions, and to degenerate mobilities. Its central thrust is the upwind methodology, which we combine with a semi-implicit formulation for the free-energy terms based on the classical convex-splitting approach. The extension of the schemes to an arbitrary number of dimensions is straightforward thanks to their dimensionally split nature, which allows to efficiently solve higher-dimensional problems with a simple parallelisation. The numerical schemes are validated and tested through a variety of examples, in different dimensions, and with various contact angles between droplets and substrates.
\end{abstract}
\subjectclassification{\subjectPDF}
\tableofcontents
\listofreviews

\begin{acronym}
	\acro{CH}[CH]{Cahn-Hilliard}
\end{acronym} \section{Introduction}\label{sec:intro}

The \ac{CH} equation is a popular phase-field model initially proposed in \cite{CH1958} to describe the process of phase separation in binary alloys. Since then, it has found innumerable applications, from capillarity--wetting phenomena \cite{AVP2019,SNK2013} and diblock copolymer molecules \cite{ZCY2020} to tumour growth \cite{WVV2014,GLN2018}, image inpainting \cite{BEG2007,BHS2009,CKL2021} and topology optimization \cite{ZW2007}; see the review \cite{KLC2016}.

Like all phase-field models, the \ac{CH} equation avoids the explicit treatment of sharp interfaces altogether via thin transition regions through which pertinent variables and physical properties vary rapidly but continuously. It has a gradient-flow structure of the form
\begin{equation}\label{eq:ch}
	\pder{\phi}{t} = \div\prt*{ M(\phi) \grad \vder{\F[\phi]}{\phi}
	}
\end{equation}
where $\phi$ is the \textit{phase-field}, a continuous function of time and space which plays the role of an order parameter describing the phases of the system. In a binary system, the limiting values $\phi=1$ and $\phi=-1$ represent each of the two phases. The \textit{mobility} $M(\phi)$, is usually \textit{degenerate} with zeros at $\phi=\pm 1$,
\begin{equation}\label{eq:mobil}
	M(\phi)=M_0 \prt{1-\phi}\prt{1+\phi},
\end{equation}
but it can also be taken as a constant, $M(\phi)=M_0$~\cite{AVP2019}.

\begin{figure}
	\centering
	\includegraphics{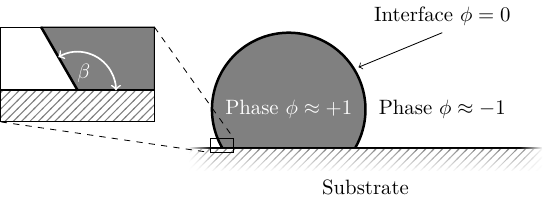}
	\caption{
		Diagram of a sessile droplet with $\phi\approx 1$ on a solid substrate and surrounded by a fluid with $\phi\approx -1$. The contact angle between the droplet and substrate is defined as $\beta$.}
	\label{fig:wall}
\end{figure}

The \textit{free energy} $\F[\phi]$ of the solution to \cref{eq:ch} is given by
\begin{equation}\label{eq:freeenergy}
	\F[\phi]=\int_{\Omega}\prt*{
		H(\phi) + \frac{\eps^2}{2} \abs*{\nabla\phi}^2
	}\dd\Omega
	+\int_{\partial \Omega}
	f_w(\phi,\beta)
	\ds,
\end{equation}
where $H(\phi)$ is a double-well potential with minima at $\phi=\pm 1$ which correspond to the stable phases in the system, $\eps$ is a positive parameter related to the width of the diffuse interface (see, for instance, \cite{CDS2005}), and $f_w(\phi,\beta)$ is the wall free energy, a function of the phase field at the boundary parametrised by the equilibrium contact angle $\beta$ \cite{WPH2012}; see \cref{fig:wall} for a schematic of a droplet on a solid substrate with contact angle $\beta$. The variation of the free energy with respect to the phase field, $\vder{\F[\phi]}{\phi}$, is known as the \textit{chemical potential}, denoted $\xi$. The boundary conditions for \eqref{eq:ch} are a combination of the natural boundary condition for the wall free energy and the no-flux condition for the chemical potential \cite{LK2011,AVP2019},
\begin{equation}\label{eq:noflux}
	\eps^2 {\nabla} \phi \cdot n = -f_w'(\phi,\beta),\quad M(\phi) {\nabla} \xi \cdot {n} = 0,
\end{equation}
where $n$ is an inward-pointing unit vector normal to the wall and $f_w'(\phi,\beta)$ denotes the derivative of $f_w(\phi,\beta)$ with respect to the phase-field.

The form of the term $f_w(\phi,\beta)$ has received considerable attention in the literature. Early contributions considered a linear form, see \cite{Seppecher1996} for instance. Here we shall assume that the function $f_w$ has bounded second derivative on $[-1,1]$, so that it can be split into a convex part and a concave part satisfying
\begin{equation}\label{eq:wallfreeenergy}
	f_w(\phi,\beta)= f_{c,w}(\phi,\beta)-f_{e,w}(\phi,\beta),
\end{equation}
where $f_{c,w}$ and $f_{e,w}$ are convex functions. A good choice is a cubic polynomial \cite{SNK2013,SNS2013,SNS2013b,AVP2019}: the lowest-order polynomial which permits the minimization of the wall free energy for the bulk densities while preventing the formation of boundary layers. Our choice here is:
\begin{equation}\label{eq:wallfreeenergy_cubic}
	f_w(\phi,\beta) = \frac{\eps\sqrt{2}}{2} \cos\beta \prt*{ \frac{\phi^3}{3} - \phi }.
\end{equation}
This cubic form alleviates the discontinuity within the realm of the diffuse-interface formulation without any additional physics. We note that the convex-concave splitting of this cubic function depends on the choice of $\beta$. Here we employ
\begin{equation}
	\begin{split}
		f_{c,w}(\phi,\beta) &=
		\begin{cases}
			\frac{\eps\sqrt{2}}{2} \cos\beta \prt*{ \frac{\phi^3}{3}-\phi+\phi^2 } & \text{if } \cos\beta \geq 0, \\[2pt]
			-\frac{\eps\sqrt{2}}{2} \cos\beta \phi^2                               & \text{otherwise},            \\
		\end{cases}
	\end{split}
\end{equation}
and choose $f_{e,w}(\phi,\beta) = f_{c,w}(\phi,\beta) - f_{w}(\phi,\beta)$, so that $f_{c,w}(\phi,\beta)$ and $f_{e,w}(\phi,\beta)$ are convex in the range $\phi \in [-1,1]$. We will often write $f_{w}(\phi)$, $f_{c,w}(\phi)$, and $f_{e,w}(\phi)$ (omitting the dependency on $\beta$) for brevity.

The variation of the free energy \eqref{eq:freeenergy} with respect to $\phi$ follows from
\begin{align}
	\der{}{\gamma} \F\prt*{ \phi+\gamma\psi } \bigg\rvert_{\gamma=0}
	 & =
	\int_{\Omega} \prt*{ \psi H'(\phi) +\eps^2 \nabla\phi \cdot \nabla \psi } \dd\Omega +
	\int_{\partial \Omega} \psi f_w'(\phi) \ds
	\\
	 & =
	\int_{\Omega} \prt*{ H'(\phi)-\eps^2 \Delta \phi } \psi \dd\Omega +
	\int_{\partial \Omega} \prt*{ \eps^2 \nabla \phi \cdot n + f_w'(\phi) } \psi \ds,
\end{align}
which, together with the boundary conditions \eqref{eq:noflux}, leads to the expression for the chemical potential:
\begin{equation}\label{eq:varfreeenergy}
	\xi = \vder{\F[\phi]}{\phi} = H'(\phi)-\eps^2 \Delta \phi.
\end{equation}
The free energy \eqref{eq:freeenergy} is dissipated in time due to the gradient-flow structure of the \ac{CH} equation. This is shown (at least formally) by differentiating the free energy $\F[\phi]$ with respect to time and applying the boundary conditions in \eqref{eq:noflux}, arriving to
\begin{equation}
	\frac{d}{dt}\F[\phi] = -\int_{\Omega} M(\phi) \abs*{ \nabla \vder{\F[\phi]}{\phi} }^2 \dd\Omega \leq 0.
\end{equation}

Mimicking this free-energy dissipation at the discrete level has been the aim of many numerical works on the \ac{CH} equation; a scheme endowed with such discrete property is called \textit{energy stable}. The first unconditionally energy-stable scheme was devised by Eyre \cite{Eyre1998}, who put forward the popular convex-splitting technique by which the potential $H(\phi)$ is separated in implicit contractive (convex) and explicit expansive (concave) terms. In fact, this semi-implicit formulation has proven so far to be the only method of deriving unconditional energy-stable schemes, and one can show that fully implicit schemes for non-linear systems, such as the \ac{CH} equation, are only conditionally energy-stable, depending on the time-step. This counter-intuitive finding has been recently analysed in \cite{XLW2019}, where the authors prove that the convex-splitting scheme for the \ac{CH} model is exactly the same as the fully-implicit scheme for a different model that is a non-trivial perturbation of the original CH equation, and the gain of stability is at the expense of a possible loss of accuracy. In any case, the convex-splitting technique is a mainstream ingredient in the construction of energy-stable schemes for the \ac{CH} equation, and has been successfully applied in various discretisation strategies such as finite differences \cite{Furihata2001,GWW2016}, finite volumes \cite{CP2008}, finite elements \cite{DWW2016,WVV2014,BBG2001,BN2008}, spectral methods \cite{HLT2007} and discontinuous Galerkin schemes \cite{AKW2013} (see \cite{TG2015} for an extensive review of energy-stable schemes for the \ac{CH} equation). A recent promising strategy to design energy-stable schemes is the so-called scalar
auxiliary variable \cite{SXY2018,SXY2019}.

The convex-splitting technique has been applied for different choices of the potential $H(\phi)$ \cite{TG2015}. In this work, we will consider the Ginzburg-Landau double-well potential,
\begin{equation}\label{eq:doublewellpot}
	H\GL(\phi) = \frac{1}{4} \prt*{\phi^2 - 1}^2;
\end{equation}
and the logarithmic potential $H_{\text{log}}(\phi)$,
\begin{equation}\label{eq:logarithmicpot}
	H\FH(\phi) = \frac{\theta}{2}\brk*{
		\prt{1 + \phi} \ln \prt*{ \frac{1+\phi}{2} }
		+ \prt{1 - \phi} \ln \prt*{ \frac{1-\phi}{2} }
	} + \frac{\theta_c}{2} \prt*{ 1-\phi^2 }
	\quad\text{for }\phi \in (-1,1),
\end{equation}
where $\theta$ and $\theta_c$ are positive constants, the \textit{absolute temperature} and \textit{absolute critical temperature}, respectively, with $\theta<\theta_c$. The logarithmic potential, known as the Flory-Huggins potential in the polymer science community \cite{Doi2013}, is more physically realistic than \eqref{eq:doublewellpot} because it can be mathematically derived from regular or ideal solution theories \cite{Doi2013,LHJ2014}. It is, however, singular as $\phi$ approaches $-1$ or $1$. The usual convex-splitting for the double-well \eqref{eq:doublewellpot} and logarithmic \eqref{eq:logarithmicpot} forms are given by
\begin{align}
	H\GL(\phi) & = H_{c,\text{GL}}(\phi) - H_{e,\text{GL}}(\phi) \coloneqq \frac{\phi^4+1}{4} - \frac{\phi^2}{2},
\end{align}
and
\begin{align}
	H\FH(\phi) & = H_{c,\text{log}}(\phi) - H_{e,\text{log}}(\phi) \\ & \coloneqq
	\frac{\theta}{2}\left[(1+\phi)\ln \left(\frac{1+\phi}{2}\right)+(1-\phi)\ln\left(\frac{1-\phi}{2}\right)\right]-\left(-\frac{\theta_c}{2}(1-\phi^2)\right),
\end{align}
where the functions $H_{c,\text{GL}}$, $H_{e,\text{GL}}$, $H_{c,\text{log}}$, and $H_{e,\text{log}}$ are all convex.

Another fundamental property of the \ac{CH} equation that has received much attention, both in terms of PDE analysis and the construction of numerical schemes, is the maximum principle, or boundedness of the phase-field. On the one hand, the phase-field solution of the \ac{CH} equation with logarithmic potential satisfies $|\phi|<1$ due to the singularities of \eqref{eq:logarithmicpot}, both for degenerate and constant mobilities. This has already been proven at the PDE level in various works \cite{AW2007,MZ2004}, including for degenerate mobilities of the type \eqref{eq:mobil} in \cite{BBG1999,EG1996}. On the other hand, for the double-well potential \eqref{eq:doublewellpot}, the phase-field solution might leave the interval $(-1,1)$ in the general case \cite{CJP2014,TG2015}. There are, however, two cases when one can analytically prove a maximum principle for the double-well potential \eqref{eq:doublewellpot}: first, for degenerate mobilities of the type \eqref{eq:mobil} which vanish when $\phi=\pm1$, the phase field is bounded in $|\phi|\leq1$, as shown in \cite{EG1996}; second, for general mobilities, by truncating the potential \eqref{eq:doublewellpot} with quadratic growth at infinities, as shown in \cite{CM1995} and applied in \cite{SY2010,CMS2011}.

The construction of numerical schemes with the discrete maximum principle property has attracted considerable attention, especially in recent years. The pioneering work by Copetti and Elliot \cite{CE1992} proposed a fully implicit scheme that verifies a discrete maximum principle under a condition for $\Delta t$, which depends on $\eps$ and the critical temperature $\theta_0$. More recent works have sought to derive schemes that unconditionally satisfy the discrete maximum principle. An important contribution has been by Chen \emph{et al.} \cite{CWW2019}, where a finite-difference scheme with unconditional discrete maximum principle for the logarithmic potential case is constructed for both constant and degenerate mobilities. In another recent work \cite{FRK2020}, the authors propose a flux-limiting technique based on high-order discontinuous Galerkin schemes and which unconditionally preserves global bounds for a family of PDEs including the \ac{CH} equation. Many other studies apply the truncation of the potential \eqref{eq:doublewellpot} with quadratic growth at infinities, with the objective of forcing their simulations to satisfy the discrete maximum principle when no rigorous proof can be derived \cite{CY2019,SY2010}. In spite of these noteworthy efforts, to this date there is no scheme flexible enough to satisfy the discrete maximum principle for a general family of free-energy potentials including the double-well \eqref{eq:doublewellpot} and the logarithmic potentials \eqref{eq:logarithmicpot}, and at the same time allowing for wetting boundary conditions such as \eqref{eq:noflux}.

The main thrust of this work is precisely the construction of a finite-volume scheme that unconditionally preserves both the discrete maximum principle and dissipation of the free energy, with second-order accuracy in space, and first-order accuracy in time. The scheme preserves these two fundamental properties for general potentials, including the double-well \eqref{eq:doublewellpot} and logarithmic \eqref{eq:logarithmicpot} potentials, wetting conditions such as \eqref{eq:wallfreeenergy} as well as more general wall free energies, and degenerate mobilities of the type \eqref{eq:mobil} which vanish when $\phi = \pm1$. In contrast to previous works, the scheme is not restricted to only particular choices of the free-energy potential. In the case of the double-well potential \eqref{eq:doublewellpot}, we do not rely on truncation techniques. Furthermore, the scheme is efficiently extended to higher dimensions while retaining all its properties through dimensional splitting. The computational efficiency can be further improved with a straightforward parallelisation. The present study builds naturally from our previous works aimed at designing structure-preserving finite-volume schemes for gradient flows and hydrodynamic systems, where a general free energy containing non-local interaction potentials drives the temporal evolution towards a steady state dictated by the minimizer of such free energy \cite{CCH2015,BCH2020,CKP2020,RPD2020,BCM2020,BCH2023,CCK2021}.

In \cref{sec:1Dscheme} we construct the one-dimensional scheme, and prove its properties of conservation of mass, boundedness of the phase field, and dissipation of the free energy. In \cref{sec:2Dscheme} we describe the general dimensional splitting technique, we construct it explicitly in two dimensions, and we prove that the properties of the previous section generalise in higher dimensions. Finally, in \cref{sec:experiments} we validate the order of convergence of the schemes, and we demonstrate their properties across examples in one, two, and three dimensions. \section{One-dimensional scheme}\label{sec:1Dscheme}

We introduce the one-dimensional version of the numerical scheme. In order to approximate the \ac{CH} equation \eqref{eq:ch} in one dimension, the computational domain $\Omega = [0,L]$ is divided into $M$ cells $C\i=\brk{x\imh,x\ih}$, each of uniform size $\Dx=L/M$, so that the centres of the cells are given by $x\i=(i-1)\Dx+\Dx/2$ for $i\in{1,\ldots,M}$. We approximate the solution at the discrete times $t\n = n\Dt$, for $n\in\set{1,\cdots,N}$, and $N$ chosen suitably. At such times, $\phi\i\n$ approximate the solution in the finite-volume sense:
\begin{equation}
	\phi\i\n \simeq \frac{1}{\Dx}\int_{C\i}\phi(t\n,x) \dx.
\end{equation}

We prescribe a finite-volume approximation of the \ac{CH} equation
\eqref{eq:ch},
\begin{subequations}\label{eq:1DScheme}
\begin{align+}\label{eq:fvsemi}
\frac{\phi\i\np - \phi\i\n}{\Dt} + \frac{F\ih\np-F\imh\np}{\Dx} = 0.
\end{align+}
The boundary fluxes are approximated in an upwind fashion, following the approach of \cite{CCH2015, BCH2020}:
\begin{align+}\label{eq:fluxsemi}
F\ih\np
=
M\prt{\phi\i\np,\phi\ip\np}
\pos{u\np\ih}
+
M\prt{\phi\ip\np,\phi\i\np}
\neg{u\np\ih}
\end{align+}
where the velocity $u\np\ih$ is given by the discrete negative gradient of the chemical potential,
\begin{align+}
u\ih\np = -\frac{\xi\np\ip - \xi\i\np}{\Dx},
\end{align+}
and where the positive and negative parts of $u\np\ih$ are constructed by
$\pos{x}=\max\set{x,0}$, $\neg{x}=\min\set{x,0}$.
The discretized mobility in \eqref{eq:fluxsemi} is approximated with the implicit values of the phase-field at either side of the cell boundary. It will be chosen as
\begin{align+}\label{eq:mobildiscr}
M\prt{x,y} = M_0\pos{1+x}\pos{1-y}
\quad\text{or}\quad
M\prt{x,y} = M_0,
\end{align+}
depending on the choice of mobility taken in \eqref{eq:mobil}. The proof of the unconditional boundedness of our numerical scheme in \cref{subsec:properties1D} relies heavily on this formulation. The numerical scheme developed in \cite{AVP2019} was restricted to constant mobility.

The discretized variation of the free energy $\xi\i\np$ is formulated semi-implicitly: the expansive part of the potential, $H_e(\rho)$, is taken explicitly, whereas the contractive part of the potential, $H_c(\rho)$, and the Laplacian, $\laplace \phi$, are taken implicitly. The approximation of $\xi\i\np$ is given by
\begin{align+}
\xi\i\np = H_c'(\phi\i\np) - H_e'(\phi\i\n) - \eps^2 (\laplace\phi)\i\np + \frac{1}{\Dx}W\i\nh,
\end{align+}
where $(\Delta \phi)\i$ is the discrete 1D second-order approximation of the Laplacian, defined as
\begin{align+}
(\Delta \phi)\i\np = \frac{\phi\ip\np-2\phi\i\np+\phi\np_{i-1}}{\Dx ^2},
\end{align+}
and the wetting term $W\i\nh$ is only evaluated at the
boundaries,
\begin{align+}\label{eq:wetting1D}
W\i\nh =
\begin{cases}
	f_{c,w}'(\phi\i\np) -f_{e,w}'(\phi\i\n) & \text{if }i=1\text{ or }i=M; \\[2pt]
	0                                       & \text{otherwise.}
\end{cases}
\end{align+}

The no-flux boundary conditions \eqref{eq:noflux} are implemented by
letting the numerical flux vanish at the boundaries,
\begin{align+}
F_{\nhalf}\np = F_{M+\nhalf}\np = 0,
\end{align+}
and by computing the Laplacian at the boundary using the \textit{ghost values}
\begin{align+}\label{eq:ghost}
\phi_{0}\np \coloneqq \phi_{1}\np,\quad
\phi_{M+1}\np \coloneqq \phi_{M}.
\end{align+}
\end{subequations}

\subsection{Properties of the scheme}\label{subsec:properties1D}

The finite-volume scheme \eqref{eq:1DScheme} possesses the following properties:
\begin{enumerate}
	\item \textbf{conservation of mass}:
	      $\sum_{i} \phi\i\n = \sum_{i} \phi\i\np$;
	\item \textbf{boundedness of the phase-field} $\phi$ for mobilities of the form $M(\phi)=M_0(1+\phi)(1-\phi)$:
	      if $\abs{\phi\i\n} \leq 1$ for every $i$, then $\abs{\phi\i\np} \leq 1$ for all $i$;
	\item \textbf{dissipation of the discrete free energy}, defined as
	      \begin{equation}\label{eq:freeenergydiscr}
		      \F_{\Delta}[\phi\n]
		      \coloneqq
		      \sum_{i=1}^M
		      \brk*{ H_c(\phi\i\n) - H_e(\phi\i\n) }
		      \Dx
		      +
		      \sum_{i=1}^{M-1}
		      \frac{\eps^2}{2} \abs{\prt{\nabla\phi}\ih\n}^2
		      \Dx
		      + f_w(\phi_1\n) + f_w(\phi_M\n),
	      \end{equation}
	      where $\prt{\nabla\phi}\ih\n$ is the discrete 1D approximation of the gradient at the cell interfaces, defined by $\prt{\nabla\phi}\ih\n\coloneqq \frac{\phi\ip\n-\phi\i\n}{\Dx}$; the discrete dissipation rate is
	      \begin{equation}\label{eq:discretedecay}
		      \frac{
			      \F_{\Delta}[\phi\np] - \F_{\Delta}[\phi\n]
		      }{\Dt}
		      \leq
		      -\sum_{i=1}^{M-1}
		      \min\set*{
			      M(\phi\i\np,\phi\ip\np),
			      M(\phi\ip\np,\phi\i\np)
		      }
		      \abs{u\ih\np}^2
		      \Dx
		      \leq 0.
	      \end{equation}
\end{enumerate}

\begin{proof} Some of the ideas in these proofs are inspired by the studies in \cite{BF2012,CCH2015, BCH2020, BCH2023}.
	The conservation of mass follows from integrating (summing) the finite-volume scheme \eqref{eq:fvsemi} over the computational domain and applying the no-flux conditions:
	\begin{equation}
		\sum_{i=1}^M
		\frac{\phi\i\np-\phi\i\n}{\Dt}
		\Dx
		=
		-\sum_{i=1}^M
		\frac{F\ih\np-F\imh\np}{\Dx}
		\Dx
		=
		- F_{M+\nhalf}\np + F_{\nhalf} = 0.
	\end{equation}

	To show the boundedness of our scheme, we follow the general contradiction proof in \cite{BCH2023}. Without loss of generality, assume first that there is a group of contiguous cells
	satisfying $\phi\i\np>1$, for some values $i$. The cells belonging to such group are $\{\phi_j\np, \phi_{j+1}\np,\ldots, \phi_{k}\np\}$. The proof also applies if there are more groups or if the groups have only one cell. The next step is to sum the scheme \eqref{eq:fvsemi} over the group of cells, resulting in
	\begin{equation}\label{eq:fvsemisum}
		\sum_{i=j}^k
		\frac{\phi\i\np-\phi\i\n}{\Dt}
		\Dx
		= - F_{k+\nhalf}\np + F_{j-\nhalf}\np.
	\end{equation}

	Since we have assumed that $\phi\i\np>1$ over the summed range, it follows that the left-hand side of \eqref{eq:fvsemisum} is strictly positive. As a result, the right-hand side of
	\eqref{eq:fvsemisum} also has
	to be positive,
	\begin{align}\label{eq:fvexpanded}
		0
		<
		- F_{k+\nhalf}\np + F_{j-\nhalf}\np
		 & =
		- \pos{u\np_{k+\nhalf}} M(\phi_k\np,\phi_{k+1}\np)
		- \neg{u\np_{k+\nhalf}} M(\phi_{k+1}\np,\phi_k\np)
		\\ & \quad\,
		+ \pos{u\np_{j-\nhalf}} M(\phi_{j-1}\np,\phi_{j}\np)
		+ \neg{u\np_{j-\nhalf}} M(\phi_{j}\np,\phi_{j-1}\np).
	\end{align}
	The first and fourth terms in the right-hand side of this expression are negative since the mobility function returns a non-negative value. The second and third terms are zero since $
		\phi_{k}\np> 1$ and $\phi_{j}\np> 1$. As a result, the whole of the right-hand side of \eqref{eq:fvexpanded} is negative, yielding a contradiction.

	Proving the lower bound $\phi\i\np>-1$ is done by following an identical strategy.

	To show the energy dissipation, we first multiply the scheme \eqref{eq:fvsemi} by the chemical potential $\xi\i\np$ and sum it over all cells, yielding
	\begin{equation}
		\sum_{i=1}^M \prt{ \phi\i\np - \phi\i\n } \xi\i\np = - \frac{\Dt}{\Dx}\sum_{i=1}^M \prt{ F\ih\np - F\imh\np } \xi\i\np.
	\end{equation}
	Then, by substituting the expression for $\xi\i\np$ and rearranging, it follows that
	\begin{align}\label{eq:proofsemi1}
		\sum_{i=1}^M \prt*{ \phi\i\np - \phi\i\n } \eps^2 (\Delta \phi)\i\np
		 & =
		\frac{\Dt}{\Dx} \sum_{i=1}^M \prt*{ F\ih\np - F\imh\np } \xi\i\np
		\\&\quad
		+\sum_{i=1}^M \prt*{ \phi\i\np - \phi\i\n} \prt*{ H_c'(\phi\i\np) - H_e'(\phi\i\n) } \\
		 & \quad
		+ \prt{ \phi_1\np - \phi_1\n } \frac{W_1\nh}{\Dx}
		+ \prt{ \phi_M\np - \phi_M\n } \frac{W_M\nh}{\Dx}.
	\end{align}
	These terms can be compared to the evolution of the free energy \eqref{eq:freeenergydiscr}:
	\begin{align}\label{eq:proofsemi2}
		\frac{1}{\Dx} \prt*{
			\F_{\Delta}[\phi\np] - \F_{\Delta}[\phi\n]
		}
		 & = \sum_{i=1}^M \prt*{ H_c(\phi\i\np)-H_c(\phi\i\n) }
		- \sum_{i=1}^M \prt*{ H_e(\phi\i\np)-H_e(\phi\i\n) }
		\\&\quad + \frac{\eps^2}{2} \sum_{i=1}^{M-1} \prt*{ \abs{\prt{\nabla\phi}\ih\np}^2 - \abs{\prt{\nabla\phi}\ih\n}^2 }
		\\&\quad + \frac{1}{\Dx} \prt*{ f_w(\phi_1\np) - f_w(\phi_1\n) + f_w(\phi_M\np) - f_w(\phi_M\n) } .
	\end{align}

	The term with the gradients can be estimated. The convexity of $f(s) = \frac{1}{2}\abs{s}^2$ leads to
	\begin{align}
		\frac{\eps^2}{2}
		\sum_{i=1}^{M-1}
		\prt*{
			\abs{\prt{\nabla\phi}\ih\np}^2 - \abs{\prt{\nabla\phi}\ih\n}^2
		} & \leq
		\eps^2
		\sum_{i=1}^{M-1}
		\prt{\nabla\phi}\ih\np
		\prt*{\prt{\nabla\phi}\ih\np - \prt{\nabla\phi}\ih\n}
		\\ & =
		\eps^2
		\sum_{i=1}^{M-1}
		\frac{\phi\ip\np-\phi\i\np}{\Dx}
		\prt*{\frac{\phi\ip\np-\phi\i\np}{\Dx}-\frac{\phi\ip\n-\phi\i\n}{\Dx}}.
	\end{align}
	We note
	\begin{align}
		\eps^2
		\sum_{i=1}^{M-1}
		\frac{\phi\ip\np-\phi\i\np}{\Dx}
		\frac{\phi\ip\np-\phi\i\np}{\Dx}
		 & =
		-\eps^2
		\sum_{i=2}^{M-1}
		\prt*{
			\frac{\phi\ip\np-\phi\i\np}{\Dx^2}
			-
			\frac{\phi\i\np-\phi_{i-1}\np}{\Dx^2}
		}
		\phi\i\np
		\\ & \quad\,
		+\eps^2
		\frac{\phi_{M}\np-\phi_{M-1}\np}{\Dx^2}
		\phi_M\np
		-\eps^2
		\frac{\phi_{2}\np-\phi_{1}\np}{\Dx^2}
		\phi_1\np,
	\end{align}
	using summation by parts. Similarly,
	\begin{align}
		-\eps^2
		\sum_{i=1}^{M-1}
		\frac{\phi\ip\np-\phi\i\np}{\Dx}
		\frac{\phi\ip\n-\phi\i\n}{\Dx}
		 & =
		\eps^2
		\sum_{i=2}^{M-1}
		\prt*{
			\frac{\phi\ip\np-\phi\i\np}{\Dx^2}
			-
			\frac{\phi\i\np-\phi_{i-1}\np}{\Dx^2}
		}
		\phi\i\n
		\\ & \quad
		-\eps^2
		\frac{\phi_{M}\np-\phi_{M-1}\np}{\Dx^2}
		\phi_M\n
		+\eps^2
		\frac{\phi_{2}\np-\phi_{1}\np}{\Dx^2}
		\phi_1\n.
	\end{align}
	Therefore, collecting the terms above, we conclude that
	\begin{align}\label{eq:proofsemi3}
		\frac{\eps^2}{2}
		\sum_{i=1}^{M-1}
		\prt*{
			\abs{\prt{\nabla\phi}\ih\np}^2 - \abs{\prt{\nabla\phi}\ih\n}^2
		} & \leq
		-\eps^2\sum_{i=1}^{M} (\Delta \phi)\i\np \prt*{\phi\i\np-\phi\i\n}.
	\end{align}
	The last equality holds because of the ghost values used to achieve the no-flux boundary conditions \eqref{eq:ghost}.
	We can now connect \eqref{eq:proofsemi3} and \eqref{eq:proofsemi1} to obtain
	\begin{align}\label{eq:proofsemi4}
		\frac{\eps^2}{2}\sum_{i=1}^{M-1} \prt*{ |\prt{\nabla\phi}\i\np|^2-|\prt{\nabla\phi}\i\n|^2 }
		 & \leq
		-		\frac{\Dt}{\Dx} \sum_{i=1}^M \prt*{ F\ih\np-F\imh\np } \xi\i\np                  \\
		 & \quad
		-\sum_{i=1}^M \prt*{ \phi\i\np-\phi\i\n } \prt*{ H_c'(\phi\i\np)-H_e'(\phi\i\n) } \\
		 & \quad
		- \prt*{ \phi_1\np-\phi_1\n } \frac{W_1\nh}{\Dx}- \prt*{ \phi_M\np-\phi_M\n } \frac{W_M\nh}{\Dx}.
	\end{align}
	Then, \eqref{eq:proofsemi2} can be rewritten as
	\begin{align}
		\frac{1}{\Dx} \prt*{
			\F_{\Delta}[\phi\np] - \F_{\Delta}[\phi\n]
		}
		 & \leq
		\sum_{i=1}^{M} \brk*{ H_c(\phi\i\np) - H_c(\phi\i\n) - \prt*{ \phi\i\np - \phi\i\n } H_c'(\phi\i\np) }
		\\ & \quad
		- \sum_{i=1}^{M} \brk*{ H_e(\phi\i\np)-H_e(\phi\i\n)- \prt*{ \phi\i\np - \phi\i\n } H_e'(\phi\i\n) }
		\\ & \quad
		+ \frac{1}{\Dx} \brk*{ f_w(\phi_1\np) - f_w(\phi_1\n) - \prt*{ \phi_1\np - \phi_1\n } W_1\nh }
		\\ & \quad
		+ \frac{1}{\Dx} \brk*{ f_w(\phi_M\np) - f_w(\phi_M\n) - \prt*{ \phi_M\np - \phi_M\n } W_M\nh }
		\\ & \quad
		- \frac{\Dt}{\Dx} \sum_{i=1}^{M} \prt*{ F\ih\np - F\imh\np } \xi\i\np
		\\ & \eqqcolon I + II+ III + IV + V.
	\end{align}

	Each term will be controlled individually. Due to the convexity of both $H_c(\phi)$ and $H_e(\phi)$, it holds
	\begin{align}
		 & H_c(\phi\i\np)-H_c(\phi\i\n) \leq \prt*{ \phi\i\np-\phi\i\n } H_c'(\phi\i\np),
		\\& H_e(\phi\i\np)-H_e(\phi\i\n) \geq \prt*{ \phi\i\np-\phi\i\n } H_e'(\phi\i\n).
	\end{align}
	Therefore, $I \leq 0$ and $II \leq 0$.

	The convexity of $f_{c,w}$ and $f_{e,w}$ can also be exploited:
	\begin{align}
		 & f_{c,w}(\phi\i\np) - f_{c,w}(\phi\i\n) \geq \prt*{ \phi\i\np - \phi\i\n } f_{c,w}'(\phi\i\np),
		\\& f_{e,w}(\phi\i\np) - f_{e,w}(\phi\i\n) \leq \prt*{ \phi\i\np - \phi\i\n } f_{e,w}'(\phi\i\n).
	\end{align}
	As a result, it holds
	\begin{align}
		f_w(\phi\i\np) - f_w(\phi\i\n) \leq \prt*{ \phi\i\np - \phi\i\n } W\i\nh.
	\end{align}
	Therefore, $III \leq 0$ and $IV \leq 0$.

	To control $V$, we perform summation by parts and apply no-flux conditions:
	\begin{align}
		V
		 & =
		- \frac{\Dt}{\Dx} \sum_{i=1}^{M} \prt*{ F\ih\np - F\imh\np } \xi\i\np
		\\ &=
		\frac{\Dt}{\Dx} \sum_{i=1}^{M-1} F\ih\np \prt*{ \xi\i\np-\xi\ip\np }
		\\ &=
		-\Dt \sum_{i=1}^{M-1} F\ih\np u\ih\np
		\\ &=
		-\Dt \sum_{i=1}^{M-1} \brk*{
			M\prt{\phi\i\np,\phi\ip\np}
			\pos{u\np\ih}
			+
			M\prt{\phi\ip\np,\phi\i\np}
			\neg{u\np\ih}
		} u\ih\np \\
		 & \leq
		-\Dt \sum_{i=1}^{M-1}\min\set*{
			M(\phi\i\np,\phi\ip\np),
			M(\phi\ip\np,\phi\i\np)
		}
		\abs{u\ih\np}^2 \leq0,
	\end{align}
	and this is precisely the dissipation rate for the discrete free energy written in \eqref{eq:discretedecay}.

\end{proof}
\section{Two-dimensional dimensionally split scheme}\label{sec:2Dscheme}

We now introduce the two-dimensional version of the numerical scheme. The scheme from \cref{sec:1Dscheme} can be directly extended to higher dimensions, and the properties described in \cref{subsec:properties1D} generalise appropriately (as was done in \cite{BCH2020}). However, for the sake of computational efficiency, we will directly formulate the two-dimensional scheme in a dimensionally split form, which nevertheless retains all structural properties.

To approximate the \ac{CH} equation \eqref{eq:ch} in two dimensions, the computational domain $\Omega = [0,L_x]\times[0,L_y]$ is divided into $M_x\times M_y$ cells $C\ij=\brk{x\imh,x\ih}\times\brk{y\jmh,y\jh}$, each of uniform size $\Dx\times\Dy$, where $\Dx=L_x/M_x$ and $\Dy=L_y/M_y$, so that the centres of the cells are given by $\prt{x\i,y\j}$, where $x\i=(i-1)\Dx+\Dx/2$ for $i\in{1,\ldots,M_x}$ and $y\j=(j-1)\Dy+\Dy/2$ for $j\in{1,\ldots,M_y}$. We approximate the solution at the discrete times $t\n = n\Dt$, for $n\in\set{1,\cdots,N}$, and $N$ chosen suitably. At such times, $\phi\ij\n$ approximates the solution in the finite-volume sense:
\begin{equation}
	\phi\ij\n \simeq \frac{1}{\Dx\Dy}\iint_{C\ij}\phi(t\n,x,y)\dx\dy.
\end{equation}

To compute the dimensionally-split update, we first update the solution along the $x$-direction (the \textit{rows}), one value of $j$ at a time. This yields the \textit{half-update} $\phi\ij\nh$ from $\phi\ij\n$. Then, we update the solution along the $y$-direction (the \textit{columns}), one value of $i$ at a time. This yields the complete update $\phi\ij\np$.

\paragraph{Rows update}
We define $\phi\ij^{n,\,0} \coloneqq \phi\ij\n$. We compute $M_y$ row updates, one for each value of $j\in{1,\ldots,M_y}$. The row update $\phi\ij\nr$ is computed from $\phi\ij\nrm$, though we write $\phi\ij\nnr$ and $\phi\ij\nnrm$ below for clarity. The update is given by
\begin{subequations}\label{eq:2DSchemeRow}
\begin{align+}\label{eq:2DFVRow}
\phi\ij\nnr
& =
\phi\ij\nnrm -
\begin{cases}
	\frac{\Delta t}{\Delta x} \prt*{ F\ihj\nnr - F\imhj\nnr } & \text{if } j=r,   \\
	0                                                         & \text{otherwise},
\end{cases}
\\
F\ihj\nnr & = M\prt{\phi\ij\nnr,\phi\ipj\nnr} \pos{u\nnr\ihj} + M\prt{\phi\ipj\nnr,\phi\ij\nnr} \neg{u\nnr\ihj},
\\
u\ihj\nnr & = -\frac{\xi\nnr\ipj - \xi\ij\nnr}{\Dx},
\\
\xi\ij\nnr & = H_c'(\phi\ij\nnr) - H_e'(\phi\ij\nnrm) - \eps^2 (\laplace\phi)\ij\nnr + \frac{1}{\Dx}W^{x,\,r-\nhalf}\ij + \frac{1}{\Dy}W^{y,\,r-\nhalf}\ij,
\\
(\Delta \phi)\ij\nnr & = \frac{\phi\ipj\nnr-2\phi\ij\nnr+\phi\imj\nnr}{\Dx ^2} + \frac{\phi\ijp\nnrm-2\phi\ij\nnr+\phi\ijm\nnrm}{\Dy ^2};
\end{align+}
with wetting terms given by
\begin{align+}
W^{x,\,r-\nhalf}\ij & =
\begin{cases}
	f_{c,w}'(\phi\ij\nnr) -f_{e,w}'(\phi\ij\nnrm) & \text{if }i=1\text{ or }i=M_x, \\[2pt]
	0                                             & \text{otherwise},
\end{cases}
\\
W^{y,\,r-\nhalf}\ij & =
\begin{cases}
	f_{c,w}'(\phi\ij\nnr) -f_{e,w}'(\phi\ij\nnrm) & \text{if }j=1\text{ or }j=M_y, \\[2pt]
	0                                             & \text{otherwise};
\end{cases}
\end{align+}
and with no-flux conditions implemented by
\begin{align+}
F_{\nhalf,\,j}\nnr &= F_{M_x+\nhalf,\,j}\nnr = 0,
\\
\phi_{0,\,j}\nnr & \coloneqq \phi_{1,\,j}\nnr,
\quad
\phi_{M_x+1,\,j}\nnr \coloneqq \phi_{M_x,\,j}\nnr,
\quad
\phi_{i,\,0}\nnr \coloneqq \phi_{i,\,1}\nnr,
\quad
\phi_{i,\,M_y+1}\nnr \coloneqq \phi_{i,\,M_y}\nnr.
\end{align+}
\end{subequations}
The half-update is defined as $\phi\ij\nh\coloneqq\phi\ij^{n,\,M_x}$.

\begin{remark}[Dependency between the updates]\label{th:dependencyUpdates}
	We wish to remark that the implicit problem that is solved to compute each row update is truly a one-dimensional problem. The update of row $r$ only depends on the adjacent rows through the Laplacian, but this dependency is taken \textit{explicitly} with respect to the update. The component in the $y$-direction is taken as
	\begin{align}
		\frac{\phi\ijp\nnrm-2\phi\ij\nnr+\phi\ijm\nnrm}{\Dy ^2},
	\end{align}
	rather than
	\begin{align}
		\frac{\phi\ijp\nnr-2\phi\ij\nnr+\phi\ijm\nnr}{\Dy ^2}.
	\end{align}
	This update simply uses the ``most up-to-date'' information available outside of row $r$; this is in fact essential to show that the dimensionally split scheme dissipates a discrete version of the free energy, see \cref{subsec:properties2D}.
\end{remark}

\paragraph{Columns update}
We now define $\phi\ij^{n+\nhalf,\,0} \coloneqq \phi\ij\nh$, and compute $M_x$ column updates, one for each value of $i\in{1,\ldots,M_x}$. The column update $\phi\ij\nhc$ is computed from $\phi\ij\nhcm$, though again we write $\phi\ij\nnhc$ and $\phi\ij\nnhcm$ below for clarity. The update is given by
\begin{subequations}\label{eq:2DSchemeCol}
\begin{align+}\label{eq:2DFVCol}
\phi\ij\nnhc
& =
\phi\ij\nnhcm -
\begin{cases}
	\frac{\Delta t}{\Delta y} \prt*{ F\ijh\nnhc - F\ijmh\nnhc } & \text{if } i=c,   \\
	0                                                           & \text{otherwise},
\end{cases}
\\
F\ijh\nnhc & = M\prt{\phi\ij\nnhc,\phi\ijp\nnhc} \pos{u\nnhc\ijh} + M\prt{\phi\ijp\nnhc,\phi\ij\nnhc} \neg{u\nnhc\ijh},
\\
u\ijh\nnhc & = -\frac{\xi\nnhc\ijp - \xi\ij\nnhc}{\Dy},
\\
\xi\ij\nnhc & = H_c'(\phi\ij\nnhc) - H_e'(\phi\ij\nnhcm) - \eps^2 (\laplace\phi)\ij\nnhc + \frac{1}{\Dx}W^{x,\,c-\nhalf}\ij + \frac{1}{\Dy}W^{y,\,c-\nhalf}\ij,
\\
(\Delta \phi)\ij\nnhc & = \frac{\phi\ipj\nnhcm-2\phi\ij\nnhc+\phi\imj\nnhcm}{\Dx ^2} + \frac{\phi\ijp\nnhc-2\phi\ij\nnhc+\phi\ijm\nnhc}{\Dy ^2};
\end{align+}
with wetting terms given by
\begin{align+}
W^{x,\,c-\nhalf}\ij & =
\begin{cases}
	f_{c,w}'(\phi\ij\nnhc) -f_{e,w}'(\phi\ij\nnhcm) & \text{if }i=1\text{ or }i=M_x, \\[2pt]
	0                                               & \text{otherwise},
\end{cases}
\\
W^{y,\,c-\nhalf}\ij & =
\begin{cases}
	f_{c,w}'(\phi\ij\nnhc) -f_{e,w}'(\phi\ij\nnhcm) & \text{if }j=1\text{ or }j=M_y, \\[2pt]
	0                                               & \text{otherwise};
\end{cases}
\end{align+}
and with no-flux conditions implemented by
\begin{align+}
F_{i,\,\nhalf}\nnhc &= F_{i,\,M_y+\nhalf}\nnhc = 0,
\\
\phi_{0,\,j}\nr & \coloneqq \phi_{1,\,j}\nnhc,
\quad
\phi_{M_x+1,\,j}\nr \coloneqq \phi_{M_x,\,j}\nnhc,
\quad
\phi_{i,\,0}\nr \coloneqq \phi_{i,\,1}\nnhc,
\quad
\phi_{i,\,M_y+1}\nr \coloneqq \phi_{i,\,M_y}\nnhc.
\end{align+}
\end{subequations}
The final update is defined as $\phi\ij\np\coloneqq\phi\ij^{n+\nhalf,\,M_y}$.

Throughout, the functions $M(x,y)$, $\pos{x}$, and $\neg{x}$ are defined as in \cref{sec:1Dscheme}.

\subsection{Properties of the dimensionally split scheme}\label{subsec:properties2D}

The finite-volume scheme (\ref{eq:2DSchemeRow}--\ref{eq:2DSchemeCol}) possesses the following properties:
\begin{enumerate}
	\item \textbf{conservation of mass}:
	      $\sum_{i,\,j} \phi\ij\n = \sum_{i,\,j} \phi\ij\np$;
	\item \textbf{boundedness of the phase-field} $\phi$ for mobilities of the form $M(\phi)=M_0(1+\phi)(1-\phi)$:
	      if $\abs{\phi\ij\n} \leq 1$ for every $i$ and $j$, then $\abs{\phi\ij\np} \leq 1$ for all $i$ and $j$;
	\item \textbf{dissipation of the discrete free energy}, $\F_{\Delta}[\phi\np] - \F_{\Delta}[\phi\n] \leq 0$, where
	      \begin{align}\label{eq:freeenergydiscr2D}
		      \F_{\Delta}[\phi\n]
		       & =
		      \sum_{i,\,j=1}^{M_x,\,M_y}
		      \brk*{ H_c(\phi\ij\n) - H_e(\phi\ij\n) }
		      \Dx\Dy
		      \\ & \quad
		      +
		      \sum_{i,\,j=1}^{M_x-1,\,M_y}
		      \frac{\eps^2}{2} \prt*{\frac{\phi\ipj\n-\phi\ij\n}{\Dx}}^2
		      \Dx\Dy
		      +
		      \sum_{i,\,j=1}^{M_x,\,M_y-1}
		      \frac{\eps^2}{2} \prt*{\frac{\phi\ijp\n-\phi\ij\n}{\Dy}}^2
		      \Dx\Dy
		      \\ & \quad
		      +
		      \sum_{j=1}^{M_y} \brk*{
			      f_w(\phi_{1,\,j}\n) + f_w(\phi_{M_x,\,j}\n)
		      } \Dy
		      +
		      \sum_{i=1}^{M_x} \brk*{
			      f_w(\phi_{i,\,1}\n) + f_w(\phi_{i,\,M_y}\n)
		      } \Dx
		      .
	      \end{align}
\end{enumerate}

\begin{proof}
	The proof of these properties is achieved by applying the arguments of the proof in \cref{subsec:properties1D} to each of the row and column updates of scheme (\ref{eq:2DSchemeRow}--\ref{eq:2DSchemeCol}).

	The conservation of mass follows from the fact that summing \cref{eq:2DFVCol} yields $\sum_{i,\,j=1}^{M_x,\,M_y} \phi\ij\nr = \sum_{i,\,j=1}^{M_x,\,M_y} \phi\ij\nrm$; similarly, summing \cref{eq:2DFVRow} yields $\sum_{i,\,j=1}^{M_x,\,M_y} \phi\ij\nhc = \sum_{i,\,j=1}^{M_x,\,M_y} \phi\ij\nhcm$.

	The proof of boundedness of the phase-field is accomplished by repeating the contradiction strategy employed in the one-dimensional case, applied to each row and column.

	To prove the energy dissipation, we show that the energy is dissipated on each row and column update. I.e., we show $\F_\Delta[\phi\ij\nr] \leq \F_\Delta[\phi\ij\nrm]$ and $\F_\Delta[\phi\ij\nhc] \leq \F_\Delta[\phi\ij\nhcm]$, thereby concluding
	\begin{align}
		\F_\Delta[\phi\ij\np] \leq \F_\Delta[\phi\ij\n].
	\end{align}
	For details of this, we refer the reader to \cref{th:2DDisRow,th:2DDisCol} in the Appendix.
\end{proof}

\begin{remark}[Efficiency and parallelisation of the dimensionally split scheme versus the full scheme.]
	As mentioned at the beginning of the section, it is possible to construct a full implicit scheme in two (and higher) dimensions that fulfils the structural properties described above. However, finding the solution to an implicit problem in high dimension can prove extremely costly, and the dimensionally split approach is likely to be more efficient. The comparison is stark in problems with a full Jacobian matrix (see \cite{BCH2020}); in our case, given the banded structure of the Jacobian, the advantage will be smaller, but still relevant.

	However, the real advantage of the dimensionally split scheme arises from its parallelisation. We remark that scheme (\ref{eq:2DSchemeRow}--\ref{eq:2DSchemeCol}) can be reformulated to update the rows in any order, as long as each row is updated precisely once (similarly with the columns). Suppose we set to update first the rows with an odd index; the implicit values at these rows do not depend on each other, so the update can be trivially parallelised. The same is true for the even-indexed rows, and again for the columns. Under this approach, the problem becomes \textit{perfectly parallel}; the same cannot be said for the full scheme.
\end{remark}
\section{Numerical experiments}\label{sec:experiments}

We now demonstrate the accuracy and performance of our schemes in a wealth of test cases involving the one-dimensional scheme and the dimensionally split scheme in two and three dimensions. We begin in \cref{sec:validation} by verifying the order of convergence of the scheme numerically. We then study the effects of wetting boundary conditions in \cref{sec:droplets}. We continue to explore wetting phenomena in \cref{sec:merging}, specifically with regards to the merging of droplets. In \cref{sec:shrinking}, we reproduce spontaneous shrinking and collapse phenomena in isolated droplets. Finally, in \cref{sec:separation}, we explore phase separation in two and three dimensions.

Interactive versions of the simulations presented on \cref{sec:validation,sec:droplets,sec:merging,sec:shrinking,sec:separation}
are available online \cite{BCK2023Web}. Videos of the simulations can be found in the permanent repository \cite{BCK2023Fig}.

\subsection{Order of convergence}\label{sec:validation}

\begin{figure}
	\centering
	\includegraphics{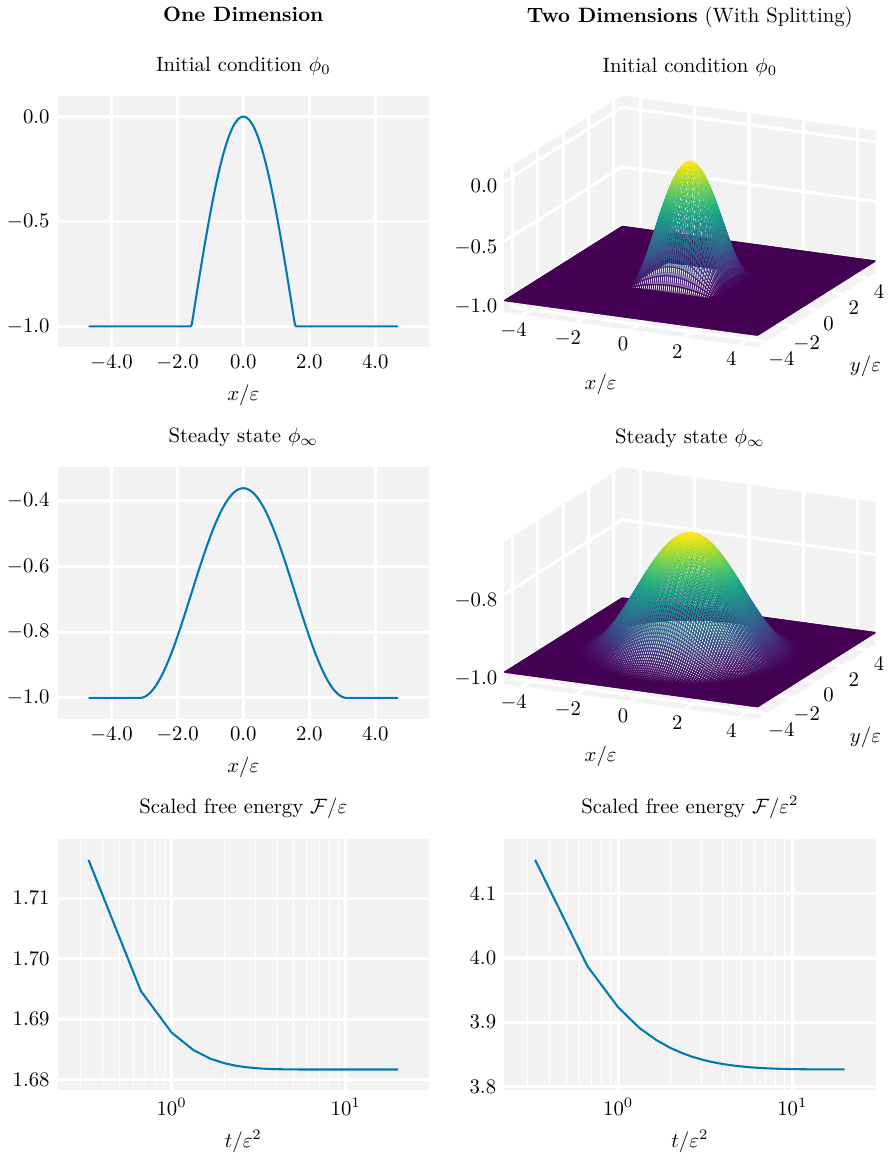}
	\caption{Initial data, stationary sates, and dissipation of the free energy for the validation test of \cref{sec:validation} in rescaled variables.
		Logarithmic potential, and degenerate mobility.
		\textbf{One dimension}: $M=256$ cells. \textbf{Two dimensions}: $M_x=M_y=128$ cells.
		See \cite{BCK2023Web,BCK2023Fig} for animations.
	}
	\label{fig:validation}
\end{figure}
\begin{figure}
	\centering
	\includegraphics{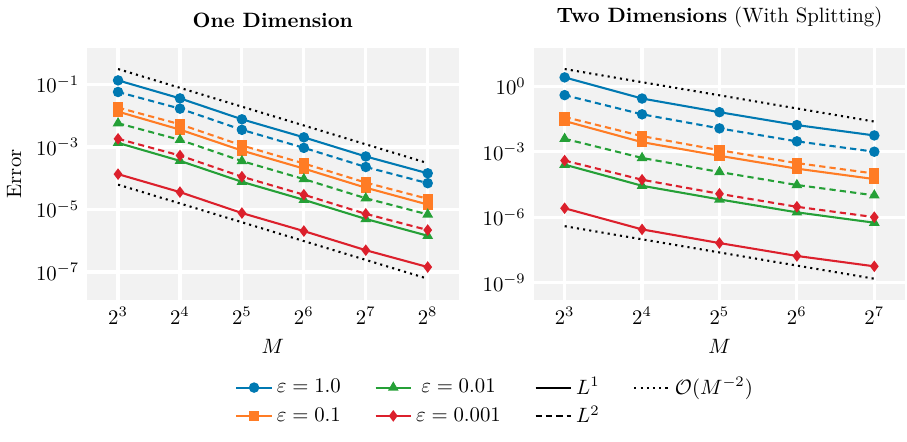}
	\caption{Order of convergence of the numerical schemes in the validation test of \cref{sec:validation}. \textbf{One dimension}: error with respect to the exact steady state \eqref{eq:steadystate1D}. \textbf{Two dimensions}: error between successive numerical solutions.
	}
	\label{fig:validationOrder}
\end{figure}

We begin the numerical experiments by numerically verifying the order of convergence of the schemes. We shall employ a test case first proposed in \cite{BBG1999} which has subsequently been used by many authors \cite{XXS2007,HHT2011}. The experiment uses the logarithmic potential \eqref{eq:logarithmicpot} in the deep quench limit (corresponding to $\theta=0$ and $\theta_c=1$, so that $H\FH(\phi)= (1-\phi^2)/2$), with the degenerate mobility \eqref{eq:mobil}. The datum for the original test case, in one dimension, is
\begin{equation}
	\phi_0(x) = \begin{cases}
		\cos\prt*{\frac{x}{\varepsilon}} - 1
		 & \text{for }
		\abs{x} \leq \frac{\pi\varepsilon}{2}
		\\
		-1
		 & \text{otherwise.}
	\end{cases}
\end{equation}
Since the \ac{CH} equation possesses symmetry and conservation of mass properties, a simple computation shows that
\begin{equation}\label{eq:steadystate1D}
	\phi_{\infty}(x) = \begin{cases}
		\frac{1}{\pi}\brk*{1+\cos\prt*{\frac{x}{\varepsilon}}} - 1
		 & \text{for }
		\abs{x} \leq \pi\varepsilon
		\\
		-1
		 & \text{otherwise}
	\end{cases}
\end{equation}
is the asymptotic steady state corresponding to $\phi_0$.

We would like to generalise this test case to two dimensions, choosing the analogous datum
\begin{equation}
	\phi_0(x,y) = \begin{cases}
		\cos\prt*{\frac{x}{\varepsilon}}
		\cos\prt*{\frac{y}{\varepsilon}}
		- 1
		 & \text{for }
		\abs{x}, \abs{y} \leq \frac{\pi\varepsilon}{2}
		\\
		-1
		 & \text{otherwise.}
	\end{cases}
\end{equation}
Annoyingly, the computation of the steady state now requires solving a two-dimensional Helmholtz equation, which can only be done semi-explicitly. For simplicity, we will perform the validation instead by computing the error between the numerical steady states of the scheme corresponding to different mesh sizes.

Crucially, we perform the validation across a wide range of values of the interface parameter: $\varepsilon=1$, $0.1$, $0.01$, and $0.001$; this ensures that the scheme performs well in different scales. In one dimension, we solve the problem on the domain $\Omega = \brk*{-\frac{3}{2} \pi\varepsilon, \frac{3}{2} \pi\varepsilon}$, and compute the solution until time $T=20\varepsilon^2$. This makes the numerical problems equivalent as $\varepsilon$ varies, up to a rescaling; in particular, the scaled evolution of the free energies (namely, $\F\varepsilon^{-1}\brk*{\phi\prt{t\varepsilon^{-2}}}$) were in agreement up to an absolute tolerance of $10^{-12}$. Similarly, we solve the two-dimensional problem on $\Omega = \brk*{-\frac{3}{2} \pi\varepsilon, \frac{3}{2} \pi\varepsilon}^2$, and on the same time interval. Again, the problems are identical up to a rescaling, and the scaled free energies ($\F\varepsilon^{-2}\brk*{\phi\prt{t\varepsilon^{-2}}}$ in this case) agreed. The initial data, stationary states, and evolution of the free energies for each problem can be found in \cref{fig:validation}.

In order to verify the second order spatial accuracy of the scheme, we choose $\Dt = \Dx^2$ throughout the test. The mesh sizes are given by $\Dx = 3\pi\varepsilon / M$, where $M=2^k$; in one dimension, $2\leq k\leq 8$; in two dimensions, $2\leq k\leq 7$. In one dimension, the error of the numerical solution is computed with respect to the exact steady state \eqref{eq:steadystate1D}; in two dimensions, the error is computed between successive numerical solutions as the mesh size is refined. The errors as a function of $M$ are shown in \cref{fig:validationOrder}. In one dimension, the scheme clearly exhibits second order accuracy; in two dimensions, the second order convergence appears to degenerate very slightly as the mesh becomes finer, likely a consequence of the dimensional splitting.
\subsection{Contact angle in wetting phenomena}\label{sec:droplets}

\begin{figure}
	\centering
	\includegraphics{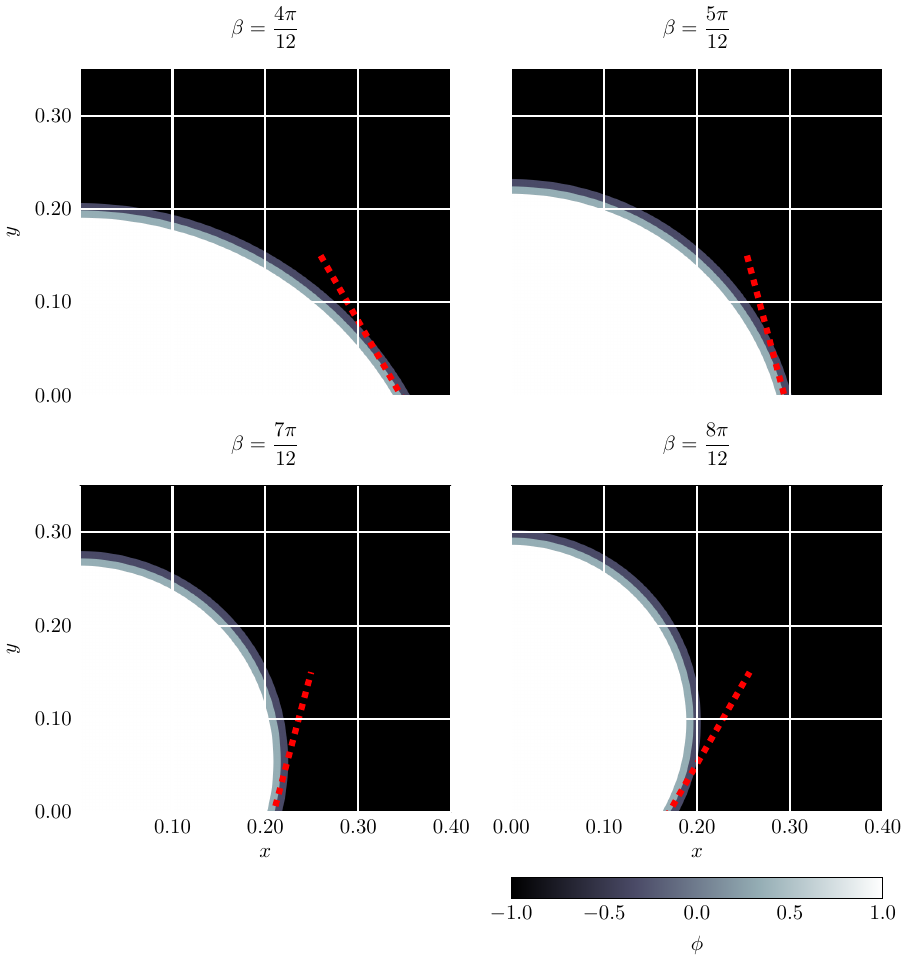}
	\caption{
		Stationary states of the ``droplet on a substrate'' tests of \cref{sec:droplets} for $\beta = 4\pi/12$, $5\pi/12$, $7\pi/12$, and $8\pi/12$.
		Double-well potential, degenerate mobility, $\varepsilon=0.01$.
		The theoretical contact angles of the droplets are shown by the dashed red lines.
		$M_x=512$, $M_y=256$, $\Dx=\Dy = 0.0015625$, $\Dt=0.01$.
		See \cite{BCK2023Web,BCK2023Fig} for animations.
	}
	\label{fig:droplets2D}
\end{figure}

We continue our journey of numerical experiments with an exploration of wetting phenomena with sessile droplets on a flat solid substrate. In this application a fluid-fluid interface moves along the solid substrate, while a contact line is formed at the intersection between the interface and the substrate. The contact angle at the three-phase conjunction, determined by the wetting properties of the substrate, is of particular interest.

In this setting, the free energy \eqref{eq:freeenergy} takes into account the wetting effects at the boundary $\partial \Omega$ between the droplet and the substrate through the function $f_w(\phi)$. Our choice of $f_w(\phi)$, \cref{eq:wallfreeenergy_cubic}, is discussed in the introduction.

We compute the stationary state of a droplet on a flat solid substrate, as depicted in \cref{fig:wall}. We solve the \ac{CH} equation on a domain $\Omega = \brk{-0.4,0.4}\times\brk{0,0.4}$ with $\varepsilon=0.01$, the double-well potential \eqref{eq:doublewellpot}, and the degenerate mobility \eqref{eq:mobil}. We initialise the droplet as a perfect semicircle,
\begin{equation}
	\phi_0(x,y) = \begin{cases}
		0.99
		 & \text{for }
		x^2+y^2 \leq 0.25^2
		\\
		-0.99
		 & \text{otherwise,}
	\end{cases}
\end{equation}
and apply the wetting boundary conditions only along the $y=0$ boundary. The fluid-fluid interface is initially perpendicular to the substrate; the wetting phenomena will drive the droplet towards a configuration where the interface meets the wall at the angle $\beta$. Similar simulations but with constant mobility have been performed in \cite{AVP2019}.

We let $M_x=512$ and $M_y=256$; this corresponds to $\Dx=\Dy = 0.0015625$. We choose $\Dt=0.01$. We solve the equation until the relative $L^1$ difference between successive iterations of the scheme is smaller than $10^{-6}$ in order to approximate the true steady state.

\Cref{fig:droplets2D} shows the contour plots of the steady states corresponding to $\beta = 4\pi/12$, $5\pi/12$, $7\pi/12$, and $8\pi/12$. The theoretical contact angle of the droplets are shown by the dashed red lines, demonstrating good agreement throughout.
\subsection{Droplet merging induced by wetting phenomena}\label{sec:merging}

\begin{figure}
	\centering
	\includegraphics{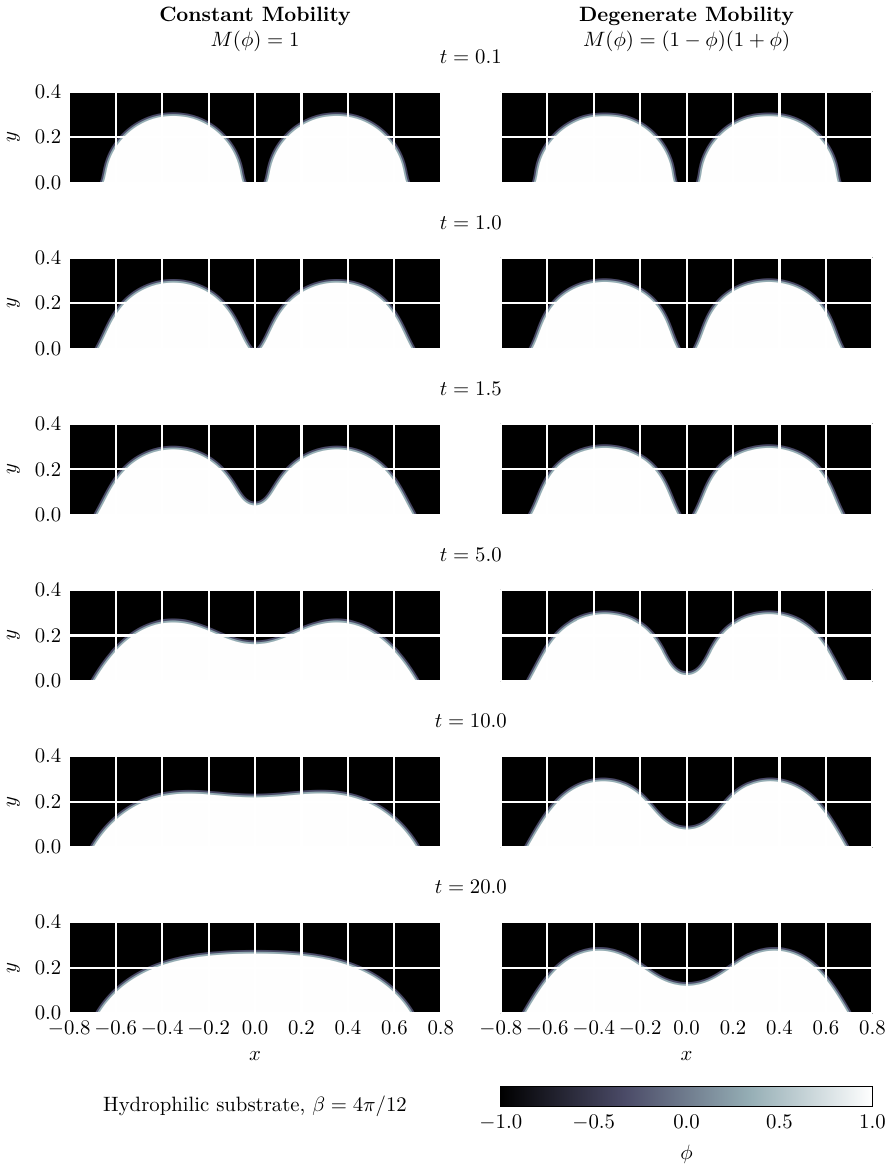}
	\caption{
		Evolution of the ``two droplets'' tests of \cref{sec:merging} for $\beta = 4\pi/12$.
		Double-well potential, $\varepsilon=0.01$.
		Left: constant mobility. Right: degenerate mobility.
		$M_x=256$, $M_y=64$, $\Dx=\Dy = 0.00625$, $\Dt=0.01$, $t\in\brk{0,20}$.
		See \cite{BCK2023Web,BCK2023Fig} for animations.
	}
	\label{fig:merging2D}
\end{figure}

\begin{figure}
	\ContinuedFloat
	\captionsetup{list=off,format=cont}
	\centering
	\includegraphics{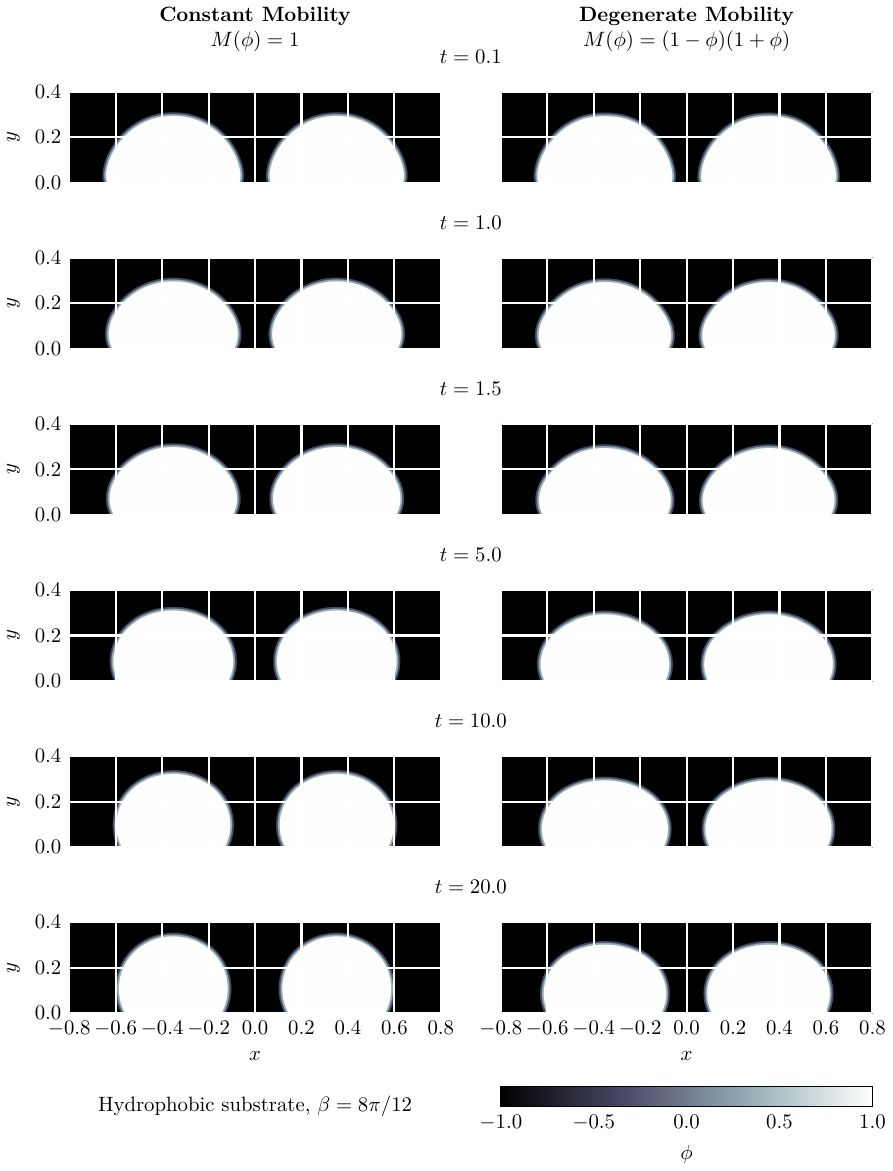}
	\caption{
		Evolution of the ``two droplets'' tests of \cref{sec:merging} for $\beta = 8\pi/12$.
		Double-well potential, $\varepsilon=0.01$.
		Left: constant mobility. Right: degenerate mobility.
		$M_x=256$, $M_y=64$, $\Dx=\Dy = 0.00625$, $\Dt=0.01$, $t\in\brk{0,20}$.
		See \cite{BCK2023Web,BCK2023Fig} for animations.
	}
\end{figure}

\begin{figure}
	\centering
	\includegraphics{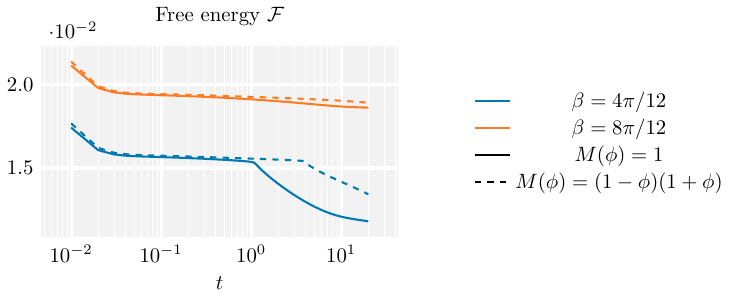}
	\caption{Evolution of the free energy in the ``two droplets'' tests of \cref{sec:merging} (see also \cref{fig:merging2D}).
		The kinks present only in the $\beta = 4\pi/12$ curves correspond to the merging of droplets.
	}
	\label{fig:merging2DEnergy}
\end{figure}
We continue the analysis of droplets and wetting phenomena from \cref{sec:droplets}. Here we consider two droplets instead, both on hydrophilic ($\beta<\pi/2$) and hydrophobic ($\beta>\pi/2$) substrates. The objective is to show that the hydrophilic substrate leads to the merging of droplets, whereas the hydrophobic substrate forces distant droplets to remain isolated.

We solve the \ac{CH} equation on a domain $\Omega = \brk{-0.8,0.8}\times\brk{0,0.4}$ with $\varepsilon=0.01$ and the double-well potential \eqref{eq:doublewellpot}, using both the constant and degenerate mobilities. We initialise two perfect droplets,
\begin{equation}
	\phi_0(x,y) = \begin{cases}
		0.99
		 & \text{for }
		(x-0.35)^2+y^2 \leq 0.3^2
		\text{ or }
		(x+0.35)^2+y^2 \leq 0.3^2
		\\
		-0.99
		 & \text{otherwise,}
	\end{cases}
\end{equation}
and again apply the wetting boundary conditions only along the $y=0$ boundary. We let $M_x=256$ and $M_y=64$; this corresponds to $\Dx=\Dy = 0.00625$. We choose $\Dt=0.01$. We solve the equation in the range $t\in\brk{0,20}$.

\Cref{fig:merging2D} shows the time evolution of the solution, both for the constant and degenerate mobilities, first on a hydrophilic substrate ($\beta = 4\pi/12$) and then on a hydrophobic setting ($\beta = 8\pi/12$). As expected, the hydrophilic experiment leads to the merging of the droplets, whereas they remain separated on the hydrophobic substrate. In both cases, the experiments with constant or degenerate mobilities exhibit qualitatively similar behaviours; the one with constant mobility leads to faster equilibration, as expected.

\Cref{fig:merging2DEnergy} shows the corresponding evolution of the free energy. In the hydrophobic setting, the value and dissipation rate of the free energy do not vary significantly between experiments with different mobilities. In the hydrophilic setting, however, the two experiments temporarily diverge, as the merging of the droplets (seen as a sharp kink on the energy curve) happens at different times.
\subsection{Spontaneous shrinking of droplets and droplet collapse}\label{sec:shrinking}

\begin{figure}
	\centering
	\includegraphics{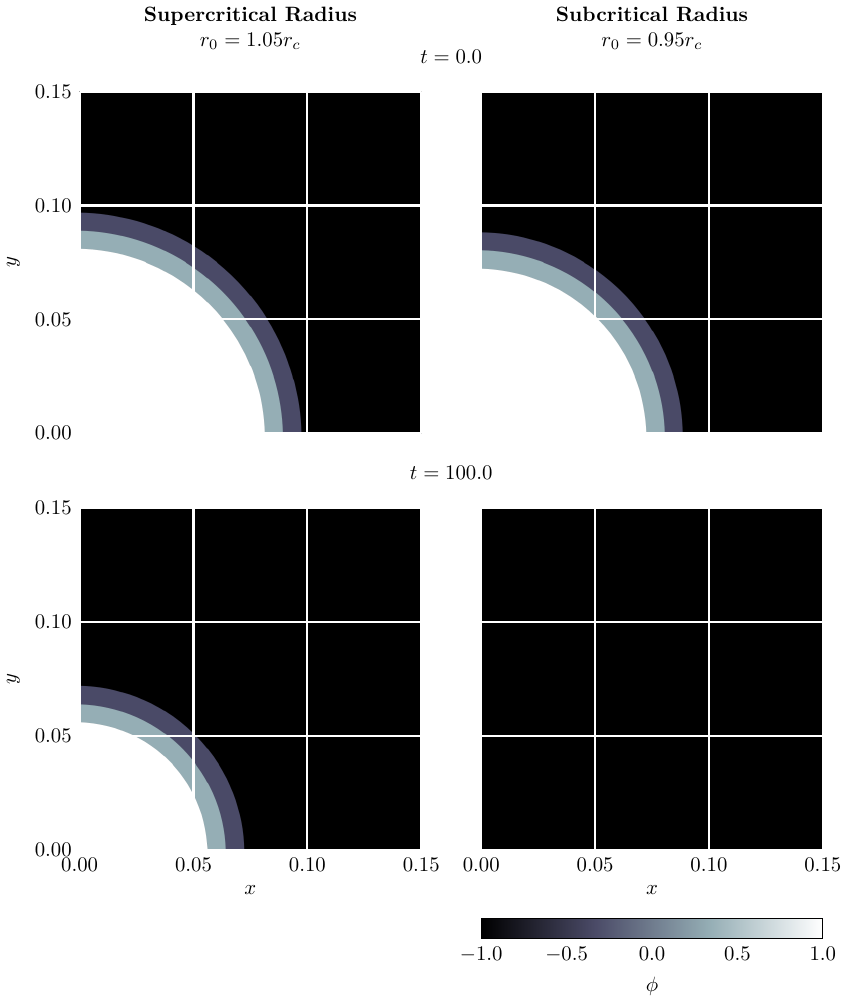}
	\caption{
		Evolution of the ``spontaneous shrinking'' tests of \cref{sec:shrinking}.
		Double-well potential, constant mobility, $\varepsilon=0.01$.
		Left: supercritical radius, $r_0=1.05r_c$. Right: subcritical radius, $r_0=0.95r_c$.
		$M_x=M_y=256$, $\Dx=\Dy = 0.00312$, $\Dt=0.01$, $t\in\brk{0,100}$.
		See \cite{BCK2023Web,BCK2023Fig} for animations.
	}
	\label{fig:shrinking}
\end{figure}

\begin{figure}
	\centering
	\includegraphics{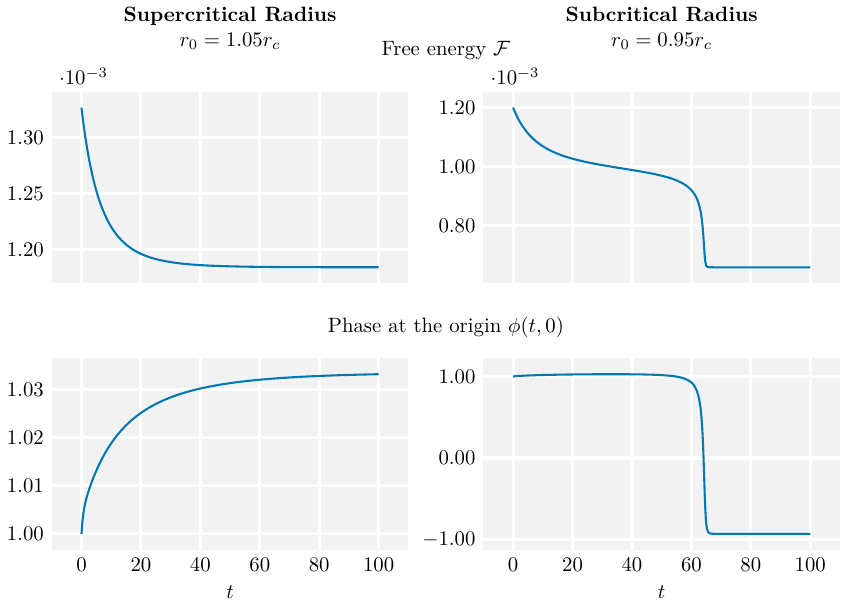}
	\caption{
		Evolution of the free energy and the phase at the origin in the ``spontaneous shrinking'' tests of \cref{sec:shrinking} (see also \cref{fig:shrinking}).
		The sharp transitions present only in the subcritical case correspond to the collapse of the droplet.
	}
	\label{fig:shrinkingEnergy}
\end{figure}

We conclude our analysis of droplet behaviour by exploring two interesting phenomena described in \cite{YZF2007}. This study considers the \ac{CH} equation with the double-well potential \eqref{eq:doublewellpot} and a constant mobility. The authors look at a single droplet centred at the origin of a square domain $\Omega=\prt{-L,L}^2$, represented by the datum
\begin{equation}\label{eq:shrinkingDatum}
	\phi_0(x,y) = \tanh\prt*{\frac{r_0 - \sqrt{x^2+y^2}}{\varepsilon\sqrt{2}}}
\end{equation}
for some $r_0>0$. They prescribe no-flux conditions along the boundary of $\Omega$.

In the regime where $\varepsilon\ll r_0$, the domain can be decomposed into the \textit{bulk} region (the inside and outside of the droplet), where $\phi \simeq \pm 1$, and the \textit{interfacial} region, where the value of $\phi$ smoothly transitions from $1$ to $-1$. The authors of \cite{YZF2007} observe that the free energy in the bulk region is negligible, since in this region the phase-field is essentially constant and takes values which minimise the potential \eqref{eq:doublewellpot}; therefore, the free energy of the datum \eqref{eq:shrinkingDatum} reduces to the contribution of the interfacial region, which is proportional to the perimeter of the droplet.

Precisely because the droplet is considered on a finite domain, it is possible to reduce its radius by shifting the value of $\phi$ in the bulk region without violating the conservation of mass. Therefore, droplets of the form \eqref{eq:shrinkingDatum} will spontaneously shrink, until they reach an equilibrium
\begin{equation}
	\phi_\infty(x,y) \simeq \tanh\prt*{\frac{r_0 - \delta r - \sqrt{x^2+y^2}}{\varepsilon\sqrt{2}}} + \delta\phi,
\end{equation}
for some $\delta r$, $\delta\phi>0$, which is more energetically favourable. Furthermore, if $r_0$ is below a critical radius $r_c$, the droplets will collapse altogether, leading to a constant equilibrium:
\begin{equation}
	\phi_\infty(x,y) = -1 + \delta\phi.
\end{equation}
In two dimensions, the formula for the critical radius is:
\begin{equation}
	r_c = \prt*{
		\frac{\sqrt{6}}{2\pi} L^2 \varepsilon
	}^{\frac{1}{3}}
\end{equation}

We solve the problem on a domain $\Omega = \brk{-0.4,0.4}^2$ with no-flux conditions and $\varepsilon=0.01$, which yields $r_c \simeq 0.086$. We consider both supercritical and subcritical cases with $r_0 = 0.95r_c$ and $r_0 = 1.05r_c$, respectively. We let $M_x=M_y=256$; this corresponds to $\Dx=\Dy = 0.00312$. We also choose $\Dt=0.01$. We compute the solution in the range $t\in\brk{0,100}$, in order to approximate the asymptotic steady state.

\Cref{fig:shrinking} displays the time evolution of the solution, for both the supercritical and subcritical cases. As expected, the supercritical droplet shrinks to a stable radius, whereas the subcritical one collapses, leading to a constant steady state. \Cref{fig:shrinkingEnergy} plots the corresponding evolution of the free energy and the value of the phase-field at the origin, $\phi(t,0)$. In the supercritical case, the free energy decays smoothly and the value of the phase-field inside the droplet tends to a value larger than $1$. On the other hand, the subcritical case initially evolves as the supercritical one, before the droplet rapidly collapses, leading to a sharp decay of the free energy before reaching the constant steady state.
\subsection{Phase separation with well-mixed initial condition}\label{sec:separation}

\begin{figure}
	\centering
	\includegraphics{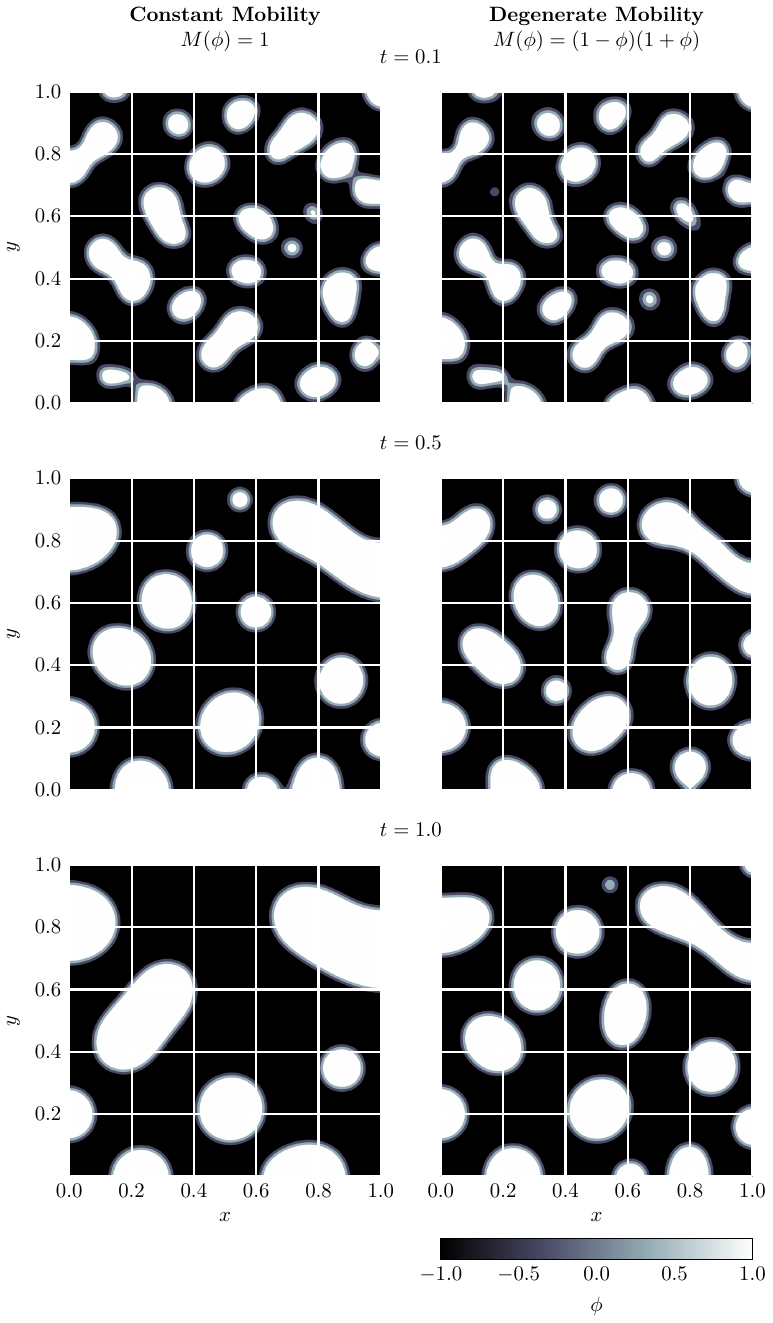}
	\caption{
		Evolution of the ``mixed initial condition'' test of \cref{sec:separation} in two dimensions.
		Double-well potential, $\varepsilon=0.01$.
		Left: constant mobility. Right: degenerate mobility.
		$M_x=M_y=128$, $\Dx=\Dy = 0.0078125$, $\Dt=0.01$, $t\in\brk{0,10}$.
		See \cite{BCK2023Web,BCK2023Fig} for animations.
	}
	\label{fig:separation2D}
\end{figure}

\begin{figure}
	\ContinuedFloat
	\captionsetup{list=off,format=cont}
	\centering
	\includegraphics{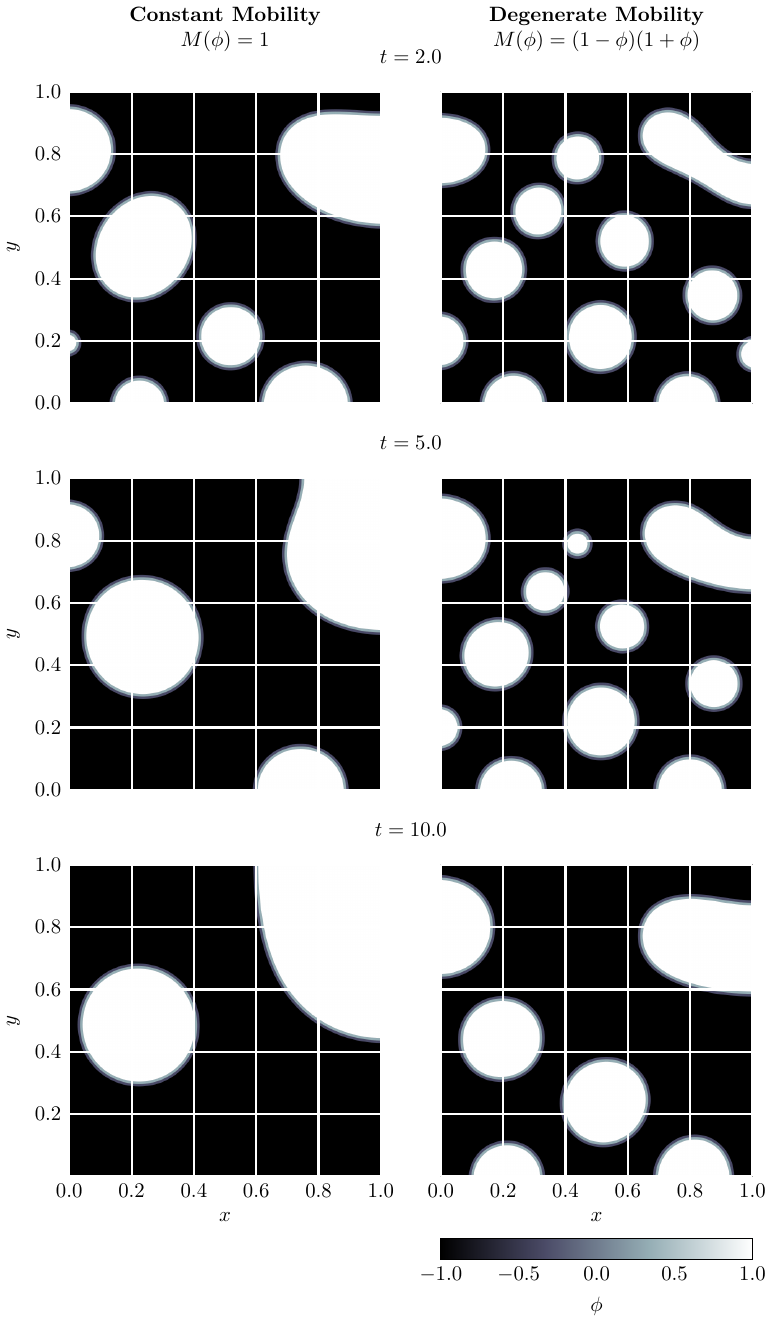}
	\caption{
		Evolution of the ``mixed initial condition'' test of \cref{sec:separation} in two dimensions.
		Double-well potential, $\varepsilon=0.01$.
		Left: constant mobility. Right: degenerate mobility.
		$M_x=M_y=128$, $\Dx=\Dy = 0.0078125$, $\Dt=0.01$, $t\in\brk{0,10}$.
		See \cite{BCK2023Web,BCK2023Fig} for animations.
	}
\end{figure} \begin{figure}
	\centering

	\begin{subfigure}[b]{0.48\textwidth}
		\centering
		\caption*{$t=1.5$\vspace{-1em}}
		\includegraphics[width=\textwidth]{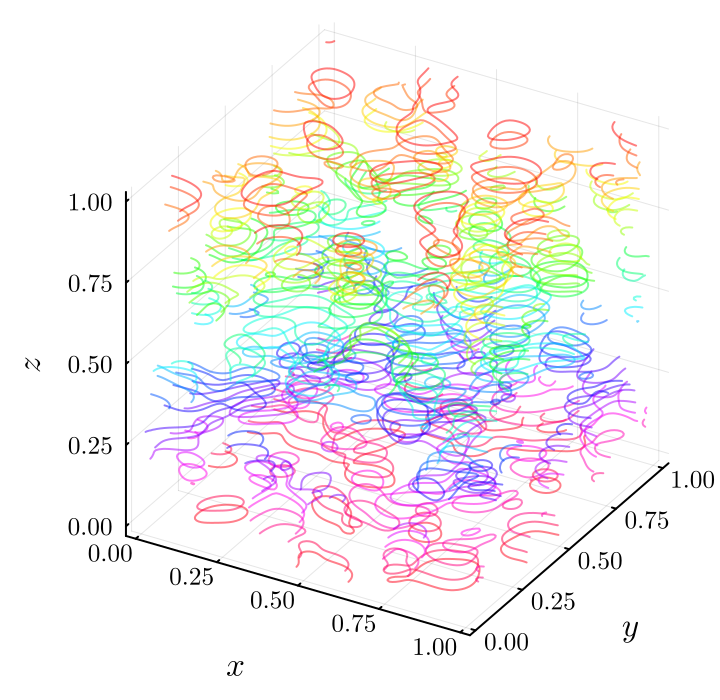}
	\end{subfigure}
	\hfill
	\begin{subfigure}[b]{0.48\textwidth}
		\centering
		\caption*{$t=5.0$\vspace{-1em}}
		\includegraphics[width=\textwidth]{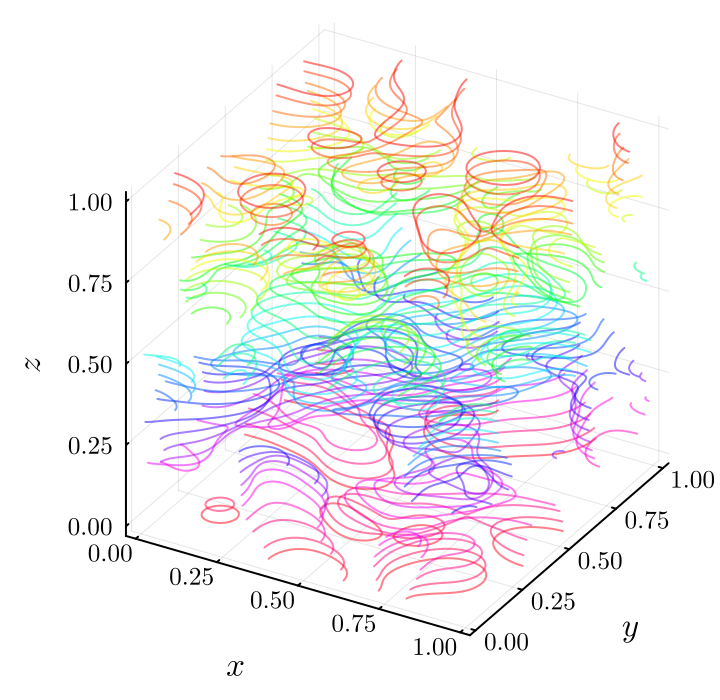}
	\end{subfigure}

	\begin{subfigure}[b]{0.48\textwidth}
		\centering
		\caption*{$t=10.0$\vspace{-1em}}
		\includegraphics[width=\textwidth]{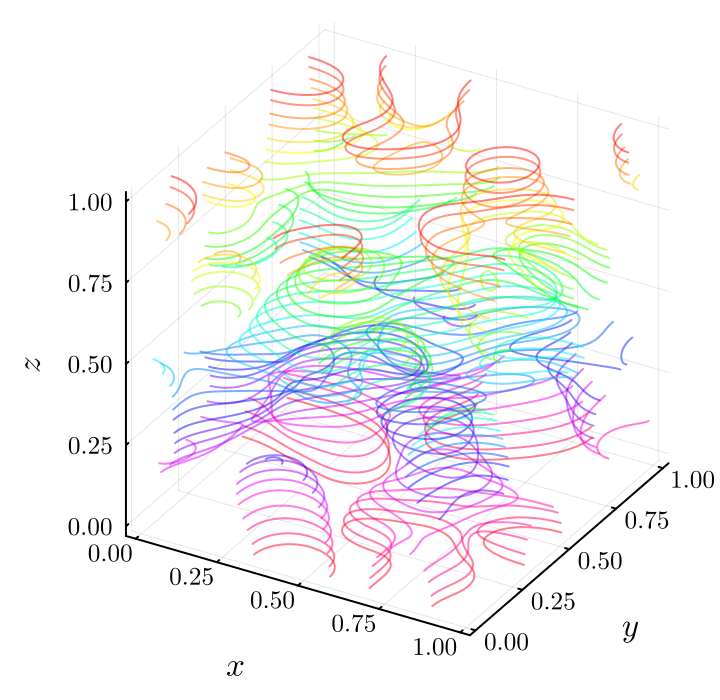}
	\end{subfigure}
	\hfill
	\begin{subfigure}[b]{0.48\textwidth}
		\centering
		\caption*{$t=20.0$\vspace{-1em}}
		\includegraphics[width=\textwidth]{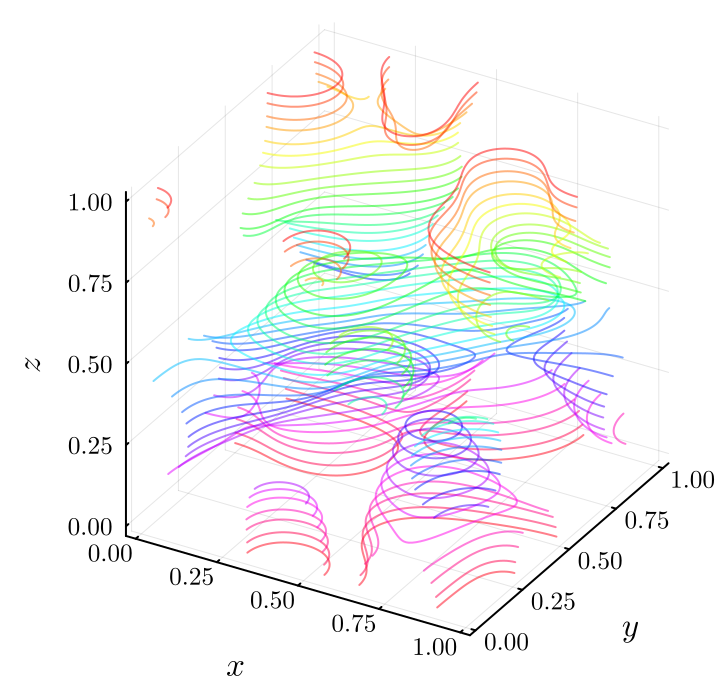}
	\end{subfigure}

	\begin{subfigure}[b]{0.48\textwidth}
		\centering
		\caption*{$t=30.0$\vspace{-1em}}
		\includegraphics[width=\textwidth]{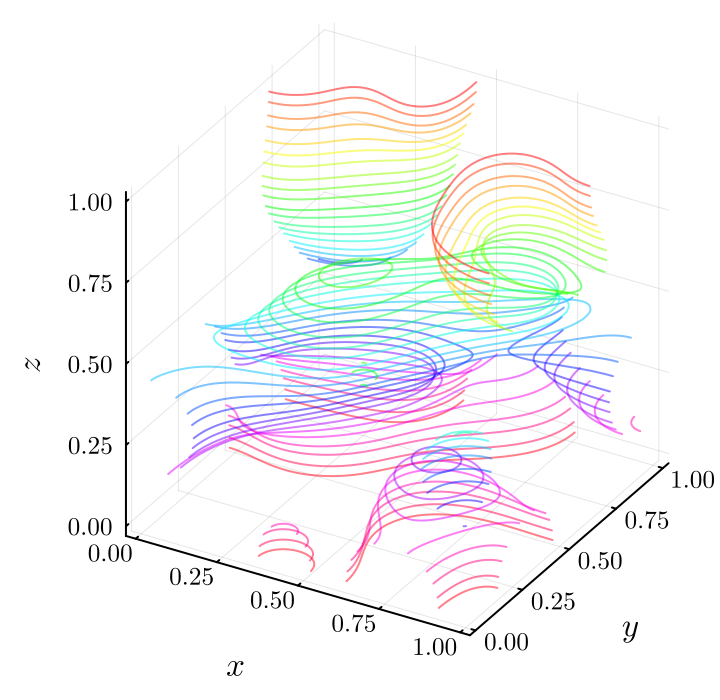}
	\end{subfigure}
	\hfill
	\begin{subfigure}[b]{0.48\textwidth}
		\centering
		\caption*{$t=40.0$\vspace{-1em}}
		\includegraphics[width=\textwidth]{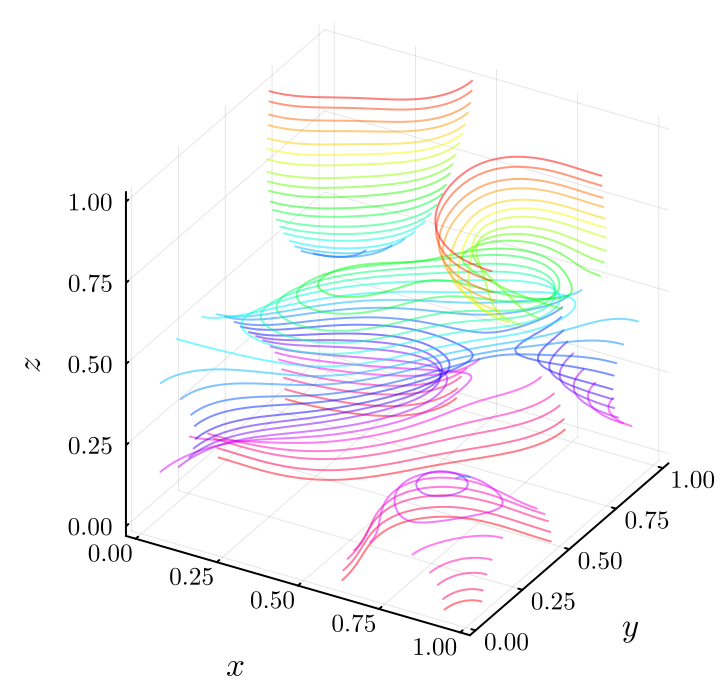}
	\end{subfigure}

	\caption{
		Evolution of the ``mixed initial condition'' test of \cref{sec:separation} in three dimensions.
		Double-well potential, constant mobility, $\varepsilon=0.01$.
		$M_x=M_y=M_z=128$, $\Dx=\Dy=\Delta z=0.0078125$, $\Dt=0.1$, $t\in\brk{0,40}$.
		See \cite{BCK2023Web,BCK2023Fig} for animations.
	}
	\label{fig:separation3D}
\end{figure} \begin{figure}
	\centering
	\includegraphics{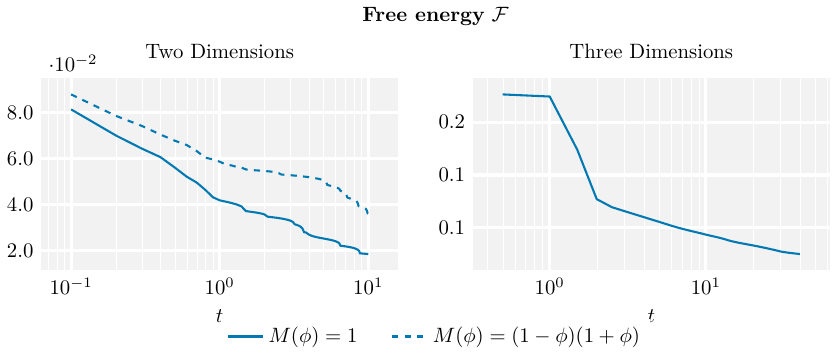}
	\caption{Evolution of the free energy in the ``mixed initial condition'' test of \cref{sec:separation} (see also \cref{fig:separation2D,fig:separation3D}).
		The kinks present on the curves correspond to the merging of droplets.
	}
	\label{fig:separationEnergy}
\end{figure}
To conclude this section, we explore the spinodal decomposition phenomenon modelled by the \ac{CH} equation. We will consider a random initial condition biased towards the $\phi=-1$ phase, in order to observe the separation of the two phases, and \textit{Ostwald ripening} (the formation of isolated droplets) in particular. On each cell of the mesh, the value of the datum $\phi_0$ will be chosen as $-0.4 + 0.25\theta$, where $\theta$ is a uniform random variable supported on $\brk{-1,1}$; this models an initial state where both phases are very well mixed. The \ac{CH} equation will then be solved with $\varepsilon=0.01$ and the double-well potential \eqref{eq:doublewellpot}.

\subsubsection{Two Dimensions}

In two dimensions, we consider a domain $\Omega=\brk{0,1}^2$, and consider both the constant and degenerate mobilities. We choose $M_x=M_y=128$, which corresponds to $\Dx=\Dy=0.0078125$, and let $\Dt = 0.01$. We solve the equation in the range $t\in\brk{0,10}$. Similar simulations can be found in \cite{BBG1999,XXS2007,GCB2008}.

\Cref{fig:separation2D} shows the time evolution of the solution, both for the constant and degenerate mobilities. As expected, small droplets first appear rapidly, and then they proceed to slowly merge into larger phases. As observed previously, the degenerate mobility significantly slows down the equilibration process.

\subsubsection{Three Dimensions}

In three dimensions, we consider a domain $\Omega=\brk{0,1}^3$, and consider just the constant mobility. We choose $M_x=M_y=M_z=128$, which corresponds to $\Dx=\Dy=\Delta z=0.0078125$, and let $\Dt = 0.1$. We solve the equation in the range $t\in\brk{0,40}$.

\Cref{fig:separation3D} shows the time evolution of the solution. Only the $\phi=0$ level sets of the solution are drawn. The colour of each contour simply indicates the value of the $z$ coordinate for ease of visualisation, and has no other physical meaning. Just as in the previous case, many small droplets arise, which then slowly merge into larger structures.

\section*{Acknowledgements}
RB was partially supported by Labex CEMPI (ANR-11-LABX-0007-01). RB and JAC were supported by the ERC Advanced Grant No. 883363 (Nonlocal PDEs for Complex Particle Dynamics (Nonlocal-CPD): Phase Transitions, Patterns and Synchronization) under the European Union’s Horizon 2020 research and innovation programme. JAC was partially supported by EPSRC Grants No. EP/V051121/1 (Stability analysis for non-linear partial differential equations across multiscale applications) under the EPSRC lead agency agreement with the NSF, and EP/T022132/1 (Spectral element methods for fractional differential equations, with applications in applied analysis and medical imaging). SK was partially supported by EPSRC Platform No. EP/L020564/1 (Multiscale Analysis of Complex Interfacial Phenomena (MACIPh): Coarse graining, Molecular modelling, stochasticity, and experimentation) and EPSRC Grant No. EP/L027186/1 (Fluid processes in smart microengineered devices: Hydrodynamics and thermodynamics in microspace). SPP acknowledges financial support from the Imperial College President’s PhD Scholarship scheme.

\renewcommand{\appendixname}{Appendix: Dissipation of the free energy in the dimensionally split scheme}

\addappendix
\label{sec:appendix}

We provide the details for the proof of energy dissipation of \cref{subsec:properties2D}.

\begin{lemma}[Dissipation, row updates]\label{th:2DDisRow}
	The finite-volume scheme (\ref{eq:2DSchemeRow}--\ref{eq:2DSchemeCol})~dissipates the discrete free energy \eqref{eq:freeenergydiscr2D} on each row update:
	\begin{align}
		\F_\Delta[\phi\ij\nr] \leq \F_\Delta[\phi\ij\nrm].
	\end{align}
\end{lemma}
\begin{proof}
	To show the energy dissipation, we first multiply the scheme \eqref{eq:2DFVRow} by the chemical potential $\xi\ij\nnr$ and sum it over all cells, yielding
	\begin{equation}
		\sum_{i,\,j=1}^{M_x,\,M_y} \prt{ \phi\ij\nnr - \phi\ij\nnrm } \xi\ij\nnr = - \frac{\Dt}{\Dx}\sum_{i=1}^{M_x} \prt{ F\ijh\nnr - F\imhj\nnr } \xi\ij\nnr.
	\end{equation}
	Then, by substituting the expression for $\xi\ij\nnr$ and rearranging, it follows that
	\begin{align}\label{eq:proofsemi1A}
		\sum_{i,\,j=1}^{M_x,\,M_y}
		\prt*{ \phi\ij\nnr - \phi\ij\nnrm } \eps^2 (\Delta \phi)\ij\nnr
		 & =
		\frac{\Dt}{\Dx} \sum_{i=1}^{M_x} \prt*{ F\ihj\nnr - F\imhj\nnr } \xi\ij\nnr
		\\&\quad
		+\sum_{i,\,j=1}^{M_x,\,M_y} \prt*{ \phi\ij\nnr - \phi\ij\nnrm} \prt*{ H_c'(\phi\ij\nnr) - H_e'(\phi\ij\nnrm) } \\
		 & \quad
		+ \sum_{i,\,j=1}^{M_x,\,M_y}
		\prt*{ \phi\ij\nnr - \phi\ij\nnrm} \brk*{
			\frac{W^{x,\,r-\nhalf}\ij}{\Dx}
			+ \frac{W^{y,\,r-\nhalf}\ij}{\Dy}
		}
	\end{align}
	These terms can be compared to the evolution of the free energy \eqref{eq:freeenergydiscr2D}:
	\begin{align}\label{eq:proofsemi2A}
		\frac{1}{\Dx\Dy} \prt*{
			\F_{\Delta}[\phi\nnr] - \F_{\Delta}[\phi\nnrm]
		}
		 & = \sum_{i,\,j=1}^{M_x,\,M_y} \prt*{ H_c(\phi\ij\nnr)-H_c(\phi\ij\nnrm) }
		- \sum_{i,\,j=1}^{M_x,\,M_y} \prt*{ H_e(\phi\ij\nnr)-H_e(\phi\ij\nnrm) }
		\\&\quad + \frac{\eps^2}{2} \sum_{i,\,j=1}^{M_x-1,\,M_y} \prt*{
			\abs*{ \frac{\phi\ipj\nnr - \phi\ij\nnr}{\Dx} }^2
			- \abs*{ \frac{\phi\ipj\nnrm - \phi\ij\nnrm}{\Dx} }^2
		}
		\\&\quad + \frac{\eps^2}{2} \sum_{i,\,j=1}^{M_x,\,M_y-1} \prt*{
			\abs*{ \frac{\phi\ijp\nnr - \phi\ij\nnr}{\Dy} }^2
			- \abs*{ \frac{\phi\ijp\nnrm - \phi\ij\nnrm}{\Dy} }^2
		}
		\\&\quad
		+
		\frac{1}{\Dx}
		\sum_{j=1}^{M_y} \brk*{
			f_w(\phi_{1,\,j}\nnr)
			- f_w(\phi_{1,\,j}\nnrm)
			+ f_w(\phi_{M_x,\,j}\nnr)
			- f_w(\phi_{M_x,\,j}\nnrm)
		}
		\\&\quad
		+
		\frac{1}{\Dy}
		\sum_{i=1}^{M_x} \brk*{
			f_w(\phi_{i,\,1}\nnr)
			- f_w(\phi_{i,\,1}\nnrm)
			+ f_w(\phi_{i,\,M_y}\nnr)
			- f_w(\phi_{i,\,M_y}\nnrm)
		}.
	\end{align}

	The term with the gradients can be estimated. We introduce the notation
	\begin{align}
		(\grad_x\phi)\ihj\nnr \coloneqq \frac{\phi\ipj\nnr - \phi\ij\nnr}{\Dx}
		\quad \text{and} \quad
		(\grad_y\phi)\ijh\nnr \coloneqq \frac{\phi\ijp\nnr - \phi\ij\nnr}{\Dx}.
	\end{align}
	We exploit the convexity of $f(s) = \frac{1}{2}\abs{s}^2$ and apply summation by parts, leading to
	\begin{align}\label{eq:proofsemi3A}
		 & \quad
		\frac{\eps^2}{2} \sum_{i,\,j=1}^{M_x-1,\,M_y} \prt*{
			\abs{ (\grad_x\phi)\ihj\nnr }^2
			- \abs{ (\grad_x\phi)\ihj\nnrm }^2
		}
		\\ & \leq
		\eps^2
		\sum_{i,\,j=1}^{M_x-1,\,M_y}
		(\grad_x\phi)\ihj\nnr
		\prt*{ (\grad_x\phi)\ihj\nnr - (\grad_x\phi)\ihj\nnrm }
		\\ & =
		\eps^2
		\sum_{i,\,j=1}^{M_x-1,\,M_y}
		\frac{\phi\ipj\nnr - \phi\ij\nnr}{\Dx}
		\prt*{\frac{\phi\ipj\nnr - \phi\ij\nnr}{\Dx} - \frac{\phi\ipj\nnrm - \phi\ij\nnrm}{\Dx}}
		\\ & =
		-\eps^2
		\sum_{i=2,\,j=1}^{M_x-1,\,M_y}
		\prt*{
			\frac{\phi\ipj\nnr-\phi\ij\nnr}{\Dx^2}
			-
			\frac{\phi\ij\nnr-\phi\imj\nnr}{\Dx^2}
		}
		\prt*{ \phi\ij\nnr-\phi\ij\nnrm }
		\\ & \quad
		+\eps^2
		\sum_{j=1}^{M_y}\brk*{
			\frac{\phi_{M,\,j}\nnr-\phi_{M-1,\,j}\nnr}{\Dx^2}
			\prt*{\phi_{M,\,j}\nnr-\phi_{M,\,j}\nnrm}
			-
			\frac{\phi_{2,\,j}\nnr-\phi_{1,\,j}\nnr}{\Dx^2}
			\prt*{\phi_{1,\,j}\nnr-\phi_{1,\,j}\nnrm}
		}
		\\ & =
		-\eps^2\sum_{i,\,j=1}^{M_x,\,M_y} \frac{\phi\ipj\nnr-2\phi\ij\nnr+\phi\imj\nnr}{\Dx ^2} \prt*{\phi\i\nnr-\phi\i\nnrm}.
	\end{align}
	The last equality holds because of the ghost values used to achieve the no-flux boundary conditions. Similarly, we find
	\begin{align}\label{eq:proofsemi3B}
		 & \quad
		\frac{\eps^2}{2} \sum_{i,\,j=1}^{M_x,\,M_y-1} \prt*{
			\abs{ (\grad_y\phi)\ijh\nnr }^2
			- \abs{ (\grad_y\phi)\ijh\nnrm }^2
		}
		\\ & \leq
		-\eps^2\sum_{i,\,j=1}^{M_x,\,M_y} \frac{\phi\ipj\nnr-2\phi\ij\nnr+\phi\imj\nnr}{\Dy ^2} \prt*{\phi\i\nnr-\phi\i\nnrm}.
	\end{align}
	All in all,
	\begin{align}\label{eq:proofsemi3C}
		 & \quad
		\frac{\eps^2}{2}
		\sum_{i,\,j=1}^{M_x-1,\,M_y} \prt*{
			\abs{ (\grad_x\phi)\ihj\nnr }^2
			- \abs{ (\grad_x\phi)\ihj\nnrm }^2
		}
		\\ &\quad
		+\frac{\eps^2}{2}
		\sum_{i,\,j=1}^{M_x,\,M_y-1} \prt*{
			\abs{ (\grad_y\phi)\ijh\nnr }^2
			- \abs{ (\grad_y\phi)\ijh\nnrm }^2
		}
		\\ &\leq
		-\eps^2\sum_{i,\,j=1}^{M_x,\,M_y} (\laplace\phi)\ij\nnr \prt*{\phi\i\nnr-\phi\i\nnrm}.
	\end{align}

	We can now connect \eqref{eq:proofsemi3C} and \eqref{eq:proofsemi1A} to obtain
	\begin{align}\label{eq:proofsemi4A}
		 & \quad
		\frac{\eps^2}{2}
		\sum_{i,\,j=1}^{M_x-1,\,M_y} \prt*{
			\abs{ (\grad_x\phi)\ihj\nnr }^2
			- \abs{ (\grad_x\phi)\ihj\nnrm }^2
		}
		\\ &\quad
		+\frac{\eps^2}{2}
		\sum_{i,\,j=1}^{M_x,\,M_y-1} \prt*{
			\abs{ (\grad_y\phi)\ijh\nnr }^2
			- \abs{ (\grad_y\phi)\ijh\nnrm }^2
		}
		\\ &\leq
		-\frac{\Dt}{\Dx} \sum_{i=1}^{M_x} \prt*{ F\ihj\nnr - F\imhj\nnr } \xi\ij\nnr                                   \\
		 & \quad
		-\sum_{i,\,j=1}^{M_x,\,M_y} \prt*{ \phi\ij\nnr - \phi\ij\nnrm} \prt*{ H_c'(\phi\ij\nnr) - H_e'(\phi\ij\nnrm) } \\
		 & \quad
		- \sum_{i,\,j=1}^{M_x,\,M_y}
		\prt*{ \phi\ij\nnr - \phi\ij\nnrm} \brk*{
			\frac{W^{x,\,r-\nhalf}\ij}{\Dx}
			+ \frac{W^{y,\,r-\nhalf}\ij}{\Dy}
		}
	\end{align}

	Then, \eqref{eq:proofsemi2A} can be rewritten as
	\begin{align}
		\frac{\F_{\Delta}[\phi\nnr] - \F_{\Delta}[\phi\nnrm]}{\Dx\Dy}
		 & \leq
		\sum_{i,\,j=1}^{M_x,\,M_y} \brk*{ H_c(\phi\ij\nnr) - H_c(\phi\ij\nnrm) - \prt*{ \phi\ij\nnr - \phi\ij\nnrm } H_c'(\phi\ij\nnr) }
		\\ & \quad
		- \sum_{i,\,j=1}^{M_x,\,M_y} \brk*{ H_e(\phi\ij\nnr) - H_e(\phi\ij\nnrm)- \prt*{ \phi\ij\nnr - \phi\ij\nnrm } H_e'(\phi\ij\nnrm) }
		\\ & \quad
		+ \frac{1}{\Dx} \sum_{j=1}^{M_y} \brk*{
			f_w(\phi_{1,\,j}\nnr)
			- f_w(\phi_{1,\,j}\nnrm)
			- \prt*{ \phi_{1,\,j}\nnr - \phi_{1,\,j}\nnrm } W_{1,\,j}^{x,\,r-\nhalf}
		}
		\\ & \quad
		+ \frac{1}{\Dx} \sum_{j=1}^{M_y} \brk*{
			f_w(\phi_{M_x,\,j}\nnr)
			- f_w(\phi_{M_x,\,j}\nnrm)
			- \prt*{ \phi_{M,\,j}\nnr - \phi_{M,\,j}\nnrm } W_{M,\,j}^{x,\,r-\nhalf}
		}
		\\ & \quad
		+ \frac{1}{\Dy} \sum_{i=1}^{M_x} \brk*{
			f_w(\phi_{i,\,1}\nnr)
			- f_w(\phi_{i,\,1}\nnrm)
			- \prt*{ \phi_{i,\,1}\nnr - \phi_{i,\,1}\nnrm } W_{i,\,1}^{y,\,r-\nhalf}
		}
		\\ & \quad
		+ \frac{1}{\Dy} \sum_{i=1}^{M_x} \brk*{
			f_w(\phi_{i,\,M_y}\nnr)
			- f_w(\phi_{i,\,M_y}\nnrm)
			- \prt*{ \phi_{i,\,M_y}\nnr - \phi_{i,\,M_y}\nnrm } W_{i,\,M_y}^{y,\,r-\nhalf}
		}
		\\ & \quad
		- \frac{\Dt}{\Dx} \sum_{i=1}^M \prt*{ F\ihj\nnr - F\imhj\nnr } \xi\ij\nnr.
	\end{align}

	Each of the sums can be controlled individually, as in the proof in \cref{subsec:properties1D}. The first six terms are controlled by convexity. The last one, applying summation by parts. This concludes the proof.
\end{proof}

\begin{lemma}[Dissipation, column updates]\label{th:2DDisCol}
	The finite-volume scheme (\ref{eq:2DSchemeRow}--\ref{eq:2DSchemeCol}) dissipates the discrete free energy \eqref{eq:freeenergydiscr2D} on each column update:
	\begin{align}
		\F_\Delta[\phi\ij\nhc] \leq \F_\Delta[\phi\ij\nhcm].
	\end{align}
\end{lemma}
\begin{proof}
	The proof is essentially identical to the proof of \cref{th:2DDisRow}.
\end{proof}
\FloatBarrier

{
	\small
	\bibliographystyle{abbrv}
	\bibliography{./BailoCarrilloKalliadasisPerez_CahnHilliardSchemes.bib}
}
\end{document}